\documentclass[a4paper, 11pt]{scrartcl}
\pdfoutput=1

\usepackage[utf8]{luainputenc}
\usepackage[english]{babel}
\usepackage{csquotes}

\usepackage[a4paper, top=0.51in, bottom=0.79in, left=0.999in, right=0.99in]{geometry}

\usepackage[plainpages=false,pdfpagelabels,hidelinks]{hyperref}

\usepackage[%
  backend=bibtex,bibencoding=ascii,
  style=authoryear-comp, dashed=false,
  firstinits=true, uniquename=init, 
  natbib=true,
  url=true,
  doi=true,
  isbn=false,
  backref=false,
  maxnames=99,
  ]{biblatex}
\usepackage{amsmath}
\usepackage{amssymb}
\usepackage{commath}
\usepackage{mathtools}

\usepackage{amsthm}
\theoremstyle{plain}
  \newtheorem{thm}{Theorem}
  \newtheorem{lem}[thm]{Lemma}

\usepackage{graphicx}
\usepackage[small]{caption}
\usepackage{subcaption}

\usepackage{booktabs}
\usepackage{rotating}

\usepackage{xparse}

\renewcommand{\vec}[1]{\underline{#1}}
\NewDocumentCommand{\mat}{mo}{%
  \IfValueTF{#2}{%
    \underline{\underline{#1}}{#2}
  }{%
    \underline{\underline{#1}}\,
  }%
}
\newcommand{\vecfnum}{\vec{f}^\mathrm{num}}
\newcommand{\fnum}{f^\mathrm{num}}

\DeclarePairedDelimiter{\diagfences}{(}{)}
\newcommand{\diag}{\operatorname{diag}\diagfences}
\newcommand{\proj}{\operatorname{proj}}

\newcommand{\scp}[2]{\left\langle{#1, #2}\right\rangle}

\renewcommand{\epsilon}{\varepsilon}
\renewcommand{\phi}{\varphi}
\newcommand{\N}{\mathbb{N}}
\newcommand{\R}{\mathbb{R}}

\usepackage{ifthen}
\newboolean{numeric_citation}
\setboolean{numeric_citation}{false}

\ifthenelse{\boolean{numeric_citation}}{
  \newcommand{\citex}[2]{#1}
}{
  \newcommand{\citex}[2]{#2}
}

\title{Enhancing stability of correction procedure via reconstruction using
       summation-by-parts operators I: Artificial dissipation}
\author{Hendrik Ranocha, Jan Glaubitz, Philipp Öffner, Thomas Sonar}
\titlehead{\includegraphics[width=0.5\textwidth]{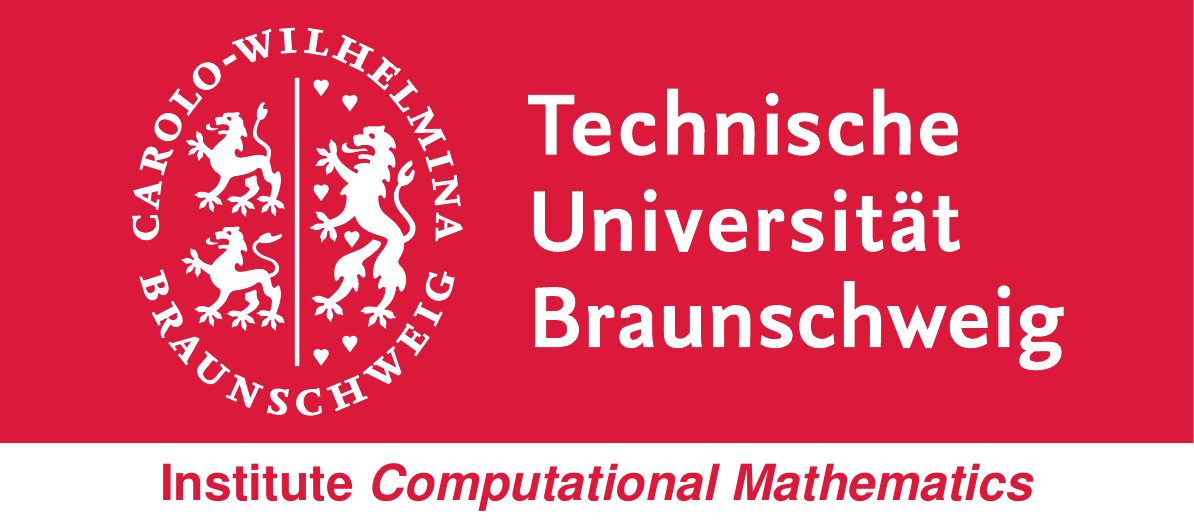}}
\date{June 2, 2016}

\hypersetup{pdfauthor={Hendrik Ranocha}}
\hypersetup{pdftitle={Artificial dissipation for CPR using SBP operators}}
\hypersetup{pdfkeywords={hyperbolic conservation laws, high order methods,
summation-by-parts, SBP, flux reconstruction, FR, lifting collocation penalty,
LCP, correction procedure via reconstruction, CPR, entropy stability,
artificial dissipation, artificial viscosity}}

\begin{document}
  \maketitle
  
  \begin{abstract}
  The correction procedure via reconstruction (CPR, also known as flux
  reconstruction) is a framework of high order semidiscretisations used for
  the numerical solution of hyperbolic conservation laws. Using a reformulation
  of these schemes relying on summation-by-parts (SBP) operators and
  simultaneous approximation terms (SATs),
  artificial dissipation / spectral viscosity operators are investigated in this
  first part of a series.
  Semidiscrete stability results for linear advection and Burgers' equation
  as model problems are extended to fully discrete stability by an explicit
  Euler method.
  As second part of this series, Glaubitz, Ranocha, Öffner, and Sonar
  (\emph{Enhancing stability of correction procedure via reconstruction using
   summation-by-parts operators II: Modal filtering}, 2016)
  investigate connections to modal filters and their application instead of
  artificial dissipation.
\end{abstract}

  \section{Introduction}

Many fundamental physical principles can be described by balance laws. Without
additional source terms, they reduce to (hyperbolic systems of) conservation laws.
These types of partial differential equations can be used inter alia as models
in fluid dynamics, electrodynamics, space and plasma physics.

Traditionally, low order numerical methods have been used to solve hyperbolic
conservation laws,
especially in industrial applications. These methods can have excellent stability
properties but also excessive numerical dissipation. Thus, they become very
costly for high accuracy or long time simulations. Therefore, in order to use
modern computing power more efficiently, high order methods provide a viable
alternative. However, these methods often lack desired stability properties.

The \emph{flux reconstruction} (FR) method has been established \citex{in}{by} \citet{huynh2007flux}
as a framework of high order semidiscretisations recovering some well known schemes
such as \emph{spectral difference} (SV) and \emph{discontinuous Galerkin} (DG)
methods with special choices of the parameters. Later, \citet{huynh2014high}
reviewed these schemes and coined the common name \emph{correction procedure via
reconstruction} (CPR).

Linearly stable schemes have been proposed \citex{in}{by} \citet{vincent2011newclass,
vincent2015extended}, extending an idea of \citet{jameson2010proof}. However,
nonlinear stability is much more difficult and has been considered inter alia
\citex{in}{by} \citet{jameson2012nonlinear, witherden2014analysis}.

\citex{Ranocha et al.}{} \citet{ranocha2016summation, ranocha2015extended, ranocha2016sbp} provided a
reformulation of CPR methods in the general framework of \emph{summation-by-parts} (SBP)
operators using \emph{simultaneous approximation terms} (SATs). These techniques
originate in \emph{finite difference} (FD) methods and have been used as building
blocks of provably stable discretisations, especially for linear (or linearised)
problems. Reviews of these schemes as well as historical and recent developments
have been published \citex{in}{by} \citet{svard2014review, nordstrom2015new, fernandez2014review}.
Generalised SBP operators have been introduced inter alia \citex{in}{by} \citet{gassner2013skew,
fernandez2014generalized} and extensions to multiple dimensions not relying on
a tensor product structure by \citet{hicken2015multidimensional, ranocha2016sbp}.

However, SBP operators with SATs have been used predominantly to create provably stable
semidiscretisations, although they can also be applied for implicit time
integration algorithms. Contrary, using a simple explicit Euler method as time
discretisation, semidiscrete stability is not sufficient for a stable fully
discrete scheme, as described inter alia in section 3.9 of \citet{ranocha2016summation}.

\emph{Artificial dissipation / spectral viscosity} has already been used in the early
works of \citex{von Neumann and Richtmyer}{} \citet{vonneumann1950method} to enhance stability of numerical schemes
for conservation laws. \citet{mattsson2004stable} used artificial dissipation
in the FD framework of SBP operators and SATs. Further developments and results
about artificial dissipation / spectral viscosity have been published inter alia
\citex{in}{by} \citet{tadmor1989convergence, ma1998chebyshev1, ma1998chebyshev2,
nordstrom2006conservative}.

In this work, the application of artificial dissipation / spectral viscosity in
the context of CPR methods using SBP operators is investigated. Therefore,
the framework of \citet{ranocha2016summation, ranocha2015extended, ranocha2016sbp}
is presented in
section \ref{sec:CPR-with-SBP}. Since SBP operators mimic integration by parts
on a discrete level, artificial dissipation operators are introduced at first in
the continuous setting in section \ref{sec:artificial-dissipation}. Using a
special representation of these viscosity operators in the semidiscrete setting,
conservation and stability are desired properties reproduced analogously.
Additionally, the discrete eigenvalues are compared with their continuous
counterparts obtained for Legendre polynomials. Fully discrete schemes are
obtained by an explicit Euler method and a new algorithm to adapt the strength of
the artificial dissipation is proposed, yielding stable fully discrete schemes
if the time step is small enough.
The theoretical results are augmented by numerical experiments in
section \ref{sec:numerical-results}. Finally, a conclusion is presented in section
\ref{sec:conclusions}, together with additional topics of further research,
inter alia the connections to modal filters described in the second part of
this series \citep{glaubitz2016enhancing}.

  \section{Correction procedure via reconstruction using summation-by-parts operators}
\label{sec:CPR-with-SBP}

In this section, the formulation of CPR methods using SBP operators given by
\citex{Ranocha, Öffner, and Sonar}{} \citet{ranocha2016summation} will be described. Additionally, $L_2$ stable
semidiscretisation will be presented for two model problems in one spatial
dimension: linear advection with constant velocity and Burgers' equation.

\subsection{Correction procedure via reconstruction}

The correction procedure via reconstruction is a semidiscretisation using a
polynomial approximation on elements. In order to describe the basic idea,
a scalar conservation law in one space dimension
\begin{equation}
\label{eq:scalar-CL}
  \partial_t u + \partial_x f(u) = 0,
\end{equation}
equipped with appropriate initial and boundary conditions will be used. Since
the focus does not lie on the implementation of boundary conditions, a compactly
supported initial condition or periodic boundary conditions will be assumed.

The domain $\Omega \subset \R$ is divided into disjoint open intervals
$\Omega_i \subset \Omega$ such that $\bigcup_i \overline\Omega_i = \Omega$.
These elements $\Omega_i$ are mapped diffeomorphically onto a standard element,
which is simply the interval $(-1,1)$ and can be realised by an affine linear
transformation in this case. In the following, all computations are performed
in this standard element.

On each element $\Omega_i$, the solution $u(t) = u(t, \cdot)$ is approximated
by a polynomial of degree $\leq p \in \N_0$. Classically, a nodal Lagrange basis
is used. Thus, the coefficients of $\vec{u}$ are nodal values
$\vec{u}_i = u(\xi_i), i \in \set{0, \dots, p}$, where $-1 \leq \xi_i \leq 1$
are distinct points in $[-1,1]$.

Since CPR methods are polynomial collocation schemes, the flux $f(u)$ is
approximated by the polynomial $\vec{f}$ interpolating at the nodes $\xi_i$,
i.e. $\vec{f}_i = f \left( \vec{u}_i \right) = f \left( u(\xi_i) \right)$.
Using a discrete derivative matrix $\mat{D}$, the divergence of $\vec{f}$ is
$\mat{D} \vec{f}$. However, since the solutions will probably be discontinuous
across elements, the discrete flux will have jump discontinuities, too.

Additionally, in order to include the effects of neighbouring elements, a
correction of the discontinuous flux is performed. Therefore, the solution
polynomial $\vec{u}$ is interpolated to the left and right boundary yielding
the values $u_L, u_R$. Thus, at each boundary node, there are two values
$u_-, u_+$ of the elements on the left and right hand side, respectively.
A continuous corrected flux with the common value of a numerical flux
$\fnum(u_-,u_+)$ (Riemann solver) at the boundary shall be used. Therefore, the
flux $\vec{f}$ is interpolated to the boundaries yielding $f_L, f_R$, and left 
and right correction functions $g_L, g_R$ are introduced. These functions are
symmetric, i.e. $g_L(x) = g_R(-x)$, approximate zero in $(-1,1)$ and fulfil
$g_L(-1) = g_R(1) = 1, g_L(1) = g_R(-1) = 0$. Finally, the time derivative
is computed via
\begin{equation}
\label{eq:CPR}
  \partial_t \vec{u}
  = - \mat{D} \vec{f}
    - (\fnum_L - f_L) \vec{g'}_L
    - (\fnum_R - f_R) \vec{g'}_R.
\end{equation}
Using a restriction matrix $\mat{R}$ performing interpolation to the boundary,
a correction matrix $\mat{C} = \left(\vec{g'}_L, \vec{g'}_R \right)$,
and writing $\vecfnum = (\fnum_L, \fnum_R)^T$, \citet{ranocha2016summation}
reformulated this as
\begin{equation}
\label{eq:SBP-CPR}
  \partial_t \vec{u}
  = - \mat{D} \vec{f}
    - \mat{C} \left( \vecfnum - \mat{R} \vec{f} \right).
\end{equation}
Using a nodal basis associated with a quadrature rule using weights
$\omega_0, \dots, \omega_p$, the mass matrix $\mat{M} = \diag{\omega_0, \dots, \omega_p}$
is diagonal and corresponds to an SBP operator, i.e. using the boundary matrix
$\mat{B} = \diag{-1,1}$, integration by parts
\begin{equation}
\label{eq:integration-by-parts}
  \int_\Omega u \, \partial_x v + \int_\Omega \partial_x u \, v
  = u \, v \big|_{\partial \Omega}
\end{equation}
is mimicked on a discrete level by summation-by-parts
\begin{equation}
  \vec{u}^T \mat{M} \mat{D} \vec{v} + \vec{u}^T \mat{D}[^T] \mat{M} \vec{v}
  = \vec{u}^T \mat{R}[^T] \mat{B} \mat{R} \vec{v},
\end{equation}
since the SBP property
\begin{equation}
\label{eq:SBP}
  \mat{M} \mat{D} + \mat{D}[^T] \mat{M}
  = \mat{R}[^T] \mat{B} \mat{R}
\end{equation}
is fulfilled. Here, either a Gauß-Legendre basis without boundary nodes
or a Lobatto-Legendre basis including both boundary nodes is used. The
associated quadrature rules are exact for polynomials of degree $\leq 2p + 1$
and $\leq 2p - 1$, respectively.

This concept has been generalised by \citet{ranocha2015extended}
to finite dimensional Hilbert spaces $X_V, X_B$ of functions on the volume
$\Omega$ and boundary $\partial \Omega$, respectively. Choosing appropriate
bases $\mathcal{B}_V$ and $\mathcal{B}_B$, the scalar product on the volume
approximates the $L_2$ scalar product
\begin{equation}
  \scp{u}{v}_M = \vec{u}^T \mat{M} \vec{v}
  \approx
  \int_\Omega u \, v = \scp{u}{v}_{L_2}
\end{equation}
and is represented by $\mat{M}$. The divergence operator is $\mat{D}$, the
scalar product on the boundary with prepended multiplication with the outer
normal is represented by $\mat{B}$, and the restriction to the boundary by
$\mat{R}$. In this way, also modal Legendre bases are allowed and fulfil the
SBP property \eqref{eq:SBP}.

The correction term $\mat{C} \left( \vecfnum - \mat{R} \vec{f} \right)$ of the
semidiscretisation \eqref{eq:SBP-CPR} corresponds to a simultaneous approximation
term in the framework of SBP methods. In general, the canonical choice of
correction matrix presented by \citet{ranocha2016summation} is
\begin{equation}
\label{eq:canonical-C}
  \mat{C} = \mat{M}[^{-1}] \mat{R}[^T] \mat{B}.
\end{equation}
However, using different choices yields the full range of energy stable schemes
derived \citex{in}{by} \citet{vincent2011newclass, vincent2015extended}.

\subsection{Linear advection}

The linear advection equation with constant velocity is a scalar conservation law
with linear flux $f(u) = u$, i.e.
\begin{equation}
\label{eq:lin-adv}
  \partial_t u + \partial_x u = 0.
\end{equation}
The canonical choice \eqref{eq:canonical-C} of the correction matrix yields the
semidiscretisation
\begin{equation}
\label{eq:lin-adv-semidisc}
  \partial_t \vec{u}
  = - \mat{D} \vec{u}
    - \mat{M}[^{-1}] \mat{R}[^T] \mat{B} \left(
        \vecfnum - \mat{R} \vec{u}
      \right)
\end{equation}
in the standard element, which is conservative across elements and stable with
respect to the discrete norm $\norm{\cdot}_M$ induced by $\mat{M}$, if an
adequate numerical flux is chosen, see inter alia
Theorem 5 of \citet{ranocha2016summation}.

\subsection{Burgers' equation}

Burgers' equation
\begin{equation}
\label{eq:Burgers}
  \partial_t u + \partial_x \frac{u^2}{2} = 0
\end{equation}
is nonlinear. Since the product of two polynomials of degree $\leq p$ is in general
a polynomial of degree $\leq 2p$, it has to be projected onto the lower dimensional
space of polynomials of degree $\leq p$. For a nodal (Gauß-Legendre or
Lobatto-Legendre) basis, the collocation approach is used, i.e. the linear
operator representing multiplication with $\vec{u}$ is given by a diagonal
matrix
\begin{equation}
  \mat{u}
  = \diag{\vec{u}}
  = \diag{\vec{u}_0, \dots, \vec{u}_p}
  = \diag{u(\xi_0), \dots, u(\xi_p)}.
\end{equation}
For a modal Legendre basis, an exact multiplication of polynomials followed by
an exact $L_2$ projection is used for the multiplication.

Using the $\mat{M}$-adjoint $\mat{u}[^*] = \mat{M}[^{-1}] \mat{u}[^T] \mat{M}$,
\citex{Ranocha et al.}{} \citet{ranocha2015extended} presented the semidiscretisation
\begin{equation}
\label{eq:Burgers-semidisc}
  \partial_t \vec{u}
  =
  - \frac{1}{3} \mat{D} \mat{u} \vec{u}
  - \frac{1}{3} \mat{u}{^*} \mat{D} \vec{u}
  + \mat{M}[^{-1}] \mat{R}[^T] \mat{B} \left(
      \vecfnum
      - \frac{1}{3} \mat{R} \mat{u} \vec{u}
      - \frac{1}{6} \left( \mat{R} \vec{u} \right)^2
    \right),
\end{equation}
which is conservative across elements and stable in the discrete norm induced
by $\mat{M}$, if an appropriate numerical flux is chosen, see inter alia
Theorem 2 of \citet{ranocha2015extended}.

  \section{Artificial dissipation / spectral viscosity}
\label{sec:artificial-dissipation}

As in the previous section, a scalar conservation law in one space dimension
\begin{equation}
  \partial_t u(t,x) + \partial_x f \left( u(t,x) \right) = 0
\end{equation}
with adequate initial and periodic boundary conditions is considered.
The introduction of a viscosity term on the right-hand side yields
\begin{equation}
\label{eq:scalar-CL-RHS}
  \partial_t u(t,x) + \partial_x f \left( u(t,x) \right)
  = (-1)^{s+1} \epsilon \left( \partial_x a(x) \partial_x \right)^{s} u(t,x),
\end{equation}
where $s \in \N$ is the \emph{order}, $\epsilon \geq 0$ the \emph{strength} and
$a \colon \R \to \R$ is a suitable function.
The term $\left( \partial_x a(x) \partial_x \right)^{s}$ describes the $s$-fold
application of the linear operator
$f(x) \mapsto \partial_x \left( a(x) \partial_x f(x) \right)$.
In the following, the dependence on $t$ and $x$ will be implied but not written
explicitly in all cases.

\subsection{Continuous setting}

In order to investigate conservation (assuming a suitably regular solution $u$),
equation \eqref{eq:scalar-CL-RHS} is integrated over some interval $\Omega$,
resulting in
\begin{align}
  \od{}{t} \int_\Omega u
  = \int_\Omega \partial_t u
  = - \int_\Omega \partial_x f(u)
    + (-1)^{s+1} \epsilon \int_\Omega \left( \partial_x a \partial_x \right)^{s} u.
\end{align}
Carrying out the integration on the right hand side yields
\begin{equation}
\begin{aligned}
  \od{}{t} \int_\Omega u
  &=- f(u) \big|_{\partial \Omega}
    + (-1)^{s+1} \epsilon \, a  \, \partial_x
      \left( \partial_x a \partial_x \right)^{s-1} u \big|_{\partial \Omega}.
\end{aligned}
\end{equation}
Thus, if $a$ vanishes at the boundary $\partial \Omega$ of the interval $\Omega$,
the rate of change of $\int_\Omega u$ is given by the flux of $f(u)$ through
the surface $\partial \Omega$, just as for the scalar conservation law
\eqref{eq:scalar-CL} with vanishing right hand side.

Investigating $L_2$ stability, equation \eqref{eq:scalar-CL-RHS} is multiplied
with the solution $u$ and integrated over $\Omega$
\begin{equation}
  \frac{1}{2} \od{}{t} \norm{u}_{L_2(\Omega)}^2
  = \frac{1}{2} \od{}{t} \int_\Omega u^2
  = \int_\Omega u \, \partial_t u
  = - \int_\Omega u \, \partial_x f(u)
    + (-1)^{s+1} \epsilon \int_\Omega u \left( \partial_x a \partial_x \right)^{s} u.
\end{equation}
Introducing the \emph{entropy flux} $F(u)$ by requiring $F'(u) = u \, f'(u)$,
the first term on the right hand side can be rewritten as
$- \int_\Omega \partial_x F(u)$, since $\partial_x F(u) = F'(u) \partial_x u
= u  f'(u) \partial_x u = u \partial_x f(u)$.
Applying integration by parts results in
\begin{equation}
\begin{aligned}
  \frac{1}{2} \od{}{t} \norm{u}_{L_2(\Omega)}^2
  =
  & - F(u) \big|_{\partial \Omega}
    + (-1)^{s+1} \epsilon \, u \, a \, \partial_x 
      \left( \partial_x a \partial_x \right)^{s-1} u \big|_{\partial \Omega}
  \\&
    + (-1)^{s} \epsilon \int_\Omega \left( a \partial_x u \right) \,
      \partial_x \left( \partial_x a \partial_x \right)^{s-1} u.
\end{aligned}
\end{equation}
Assuming again that $a$ vanishes at the boundary $\partial \Omega$, this can be
rewritten as
\begin{equation}
\begin{aligned}
  \frac{1}{2} \od{}{t} \norm{u}_{L_2(\Omega)}^2
  &=- F(u) \big|_{\partial \Omega}
    + (-1)^{s} \epsilon \int_\Omega \left( a \partial_x u \right) \,
      \partial_x \left( \partial_x a \partial_x \right)^{s-1} u
  \\
  &=- F(u) \big|_{\partial \Omega}
    + (-1)^{s+1} \epsilon \int_\Omega
      \left[ \left( \partial_x a \partial_x \right) u \right]
      \left[ \left( \partial_x a \partial_x \right)^{s-1} u \right].
\end{aligned}
\end{equation}
Using induction, this becomes
\begin{equation}
  \frac{1}{2} \od{}{t} \norm{u}_{L_2(\Omega)}^2
  = - F(u) \big|_{\partial \Omega} +
    \begin{cases}
      (-1)^{s+1} \epsilon \int_\Omega
        \left[ \left( \partial_x a \partial_x \right)^{s/2} u \right]^2
      , &s \text{ even},
      \\
      (-1)^{s} \epsilon \int_\Omega a
        \left[ \partial_x \left( \partial_x a \partial_x \right)^\frac{s-1}{2} u \right]^2
      , &s \text{ odd}.
    \end{cases}
\end{equation}
Thus, if $a$ vanishes at the boundary $\partial \Omega$, the rate of change of
the integral of the $L_2$ entropy $u \mapsto U(u) = \frac{1}{2} u^2$ is given by the 
entropy flux $F(u)$ through the surface of $\partial \Omega$ and an additional
term, which is non-positive if $a \geq 0$ in $\Omega$. Thus, the right hand side
in equation \eqref{eq:scalar-CL-RHS} has a stabilising effect.
Hence, in the spirit of a numerical method relying on an element-wise
discretisation as CPR, choosing $\Omega$ as an element and using $a \geq 0$ in
$\Omega$ with $a = 0$ on $\partial \Omega$, the right hand side of
\eqref{eq:scalar-CL-RHS} can be added as a stabilising artificial dissipation /
viscosity term not influencing conservation across elements.

However, care has to be taken during the discretisation of \eqref{eq:scalar-CL-RHS}.
Approximating the solution $u$ and the function $a$ on $\Omega$ as a polynomial
of degree $\leq p$, the exact product $a \, u$ is in general not a polynomial of
degree $\leq p$. Therefore, some kind of projection is necessary. This projection
might not be compatible with restriction of functions to the boundary, i.e. the
approximation of $(a \, u)$ might not be zero on  $\partial \Omega$ even if $a$
vanishes there. Of course, a discrete counterpart of integration by parts has
to be used: summation-by-parts.

\subsection{Semidiscrete setting}

In the following, a conservative and stable (in the discrete norm $\norm{\cdot}_M$
induced by $\mat{M}$) scheme will be augmented with an additional term, a discrete
equivalent of the viscosity term in \eqref{eq:scalar-CL-RHS}.

Using a CPR method with SBP operators, a direct discretisation of the dissipative
term, i.e. the right hand side of \eqref{eq:scalar-CL-RHS}, can be written as
\begin{equation}
\label{eq:RHS-simple}
  (-1)^{s+1} \epsilon \left( \mat{D} \mat{a} \mat{D} \right)^{s} \vec{u},
\end{equation}
where $\mat{D}$ is the derivative matrix and $\mat{a}$ represents multiplication
with $a$, followed by some projection on the space of polynomials of degree $\leq p$.
For a nodal basis (e.g. using Gauß-Legendre or Lobatto-Legendre nodes), a
collocation approach is used, i.e. $\mat{a} \vec{u}$ represents the polynomial
interpolating at the quadrature nodes. For a modal Legendre basis, an exact
$L_2$ projection will be used, as proposed \citex{in}{by} \citet{ranocha2015extended}. 

Investigating conservation across elements for $s = 1$, the semidiscrete equation
for $\partial_t \vec{u}$ is multiplied with $\vec{1}^T \mat{M}$, where
the constant function $x \mapsto 1$ is represented by $\vec{1}$. Thus, the
additional term induced by \eqref{eq:RHS-simple} divided by $\epsilon$ is ($s = 1$)
\begin{equation}
  \vec{1}^T \mat{M} \mat{D} \mat{a} \mat{D} u
  = \vec{1}^T \mat{R}[^T] \mat{B} \mat{R} \mat{a} \mat{D} u
    - \vec{1}^T \mat{D}[^T] \mat{M} \mat{a} \mat{D} u
  = \vec{1}^T \mat{R}[^T] \mat{B} \mat{R} \mat{a} \mat{D} u,
\end{equation}
where the SBP property
$\mat{M} \mat{D} = \mat{R}[^T] \mat{B} \mat{R} - \mat{D}[^T] \mat{M}$ \eqref{eq:SBP}
has been used. Since the derivative is exact for constants, $\mat{D} \vec{1} = 0$.
Thus, the resulting scheme is conservative if and only if the projection used
preserves boundary values. This is the case for a nodal Lobatto-Legendre basis
including boundary points. However, a nodal Gauß-Legendre or a modal Legendre
basis do not have this property.

Turning to stability for $s= 1$, multiplying the term \eqref{eq:RHS-simple} with
$\vec{u}^T \mat{M}$ and dividing by $\epsilon$ yields by the SBP property \eqref{eq:SBP}
\begin{equation}
    \vec{u}^T \mat{M} \mat{D} \mat{a} \mat{D} u
  = \vec{u}^T \mat{R}[^T] \mat{B} \mat{R} \mat{a} \mat{D} u
    - \vec{u}^T \mat{D}[^T] \mat{M} \mat{a} \mat{D} u.
\end{equation}
Again, the boundary term does not vanish in general. Additionally, the multiplication
matrix $\mat{a}$ has to be self-adjoint and positive semi-definite with respect
to (the scalar product induced by) $\mat{M}$ in order to guarantee that the last term
is non-positive.

However, these problems can be circumvented. Using the SBP property \eqref{eq:SBP},
the term \eqref{eq:RHS-simple} for $s = 1$ can be written as
\begin{equation}
  \epsilon \, \mat{D} \mat{a} \mat{D} u
  = \epsilon \, \mat{M}[^{-1}] \mat{M} \mat{D} \mat{a} \mat{D} u
  = \epsilon \, \mat{M}[^{-1}] \left(
      \mat{R}[^T] \mat{B} \mat{R} \mat{a} \mat{D} u
      - \mat{D}[^T] \mat{M} \mat{a} \mat{D} u
    \right).
\end{equation}
Enforcing the boundary term to vanish yields for arbitrary $s$ the discrete form
\begin{equation}
\label{eq:RHS-smart}
  (-1)^{s+1} \epsilon \left( - \mat{M}[^{-1}] \mat{D}[^T] \mat{M} \mat{a} \mat{D}
  \right)^s \vec{u}
  =
  - \epsilon \left( \mat{M}[^{-1}] \mat{D}[^T] \mat{M} \mat{a} \mat{D} \right)^{s} \vec{u}
\end{equation}
of the viscosity term. This form is similar to the one given \citex{in}{by}
\citet{mattsson2004stable} in the context of finite difference methods using SBP
operators. However, the artificial dissipation used in that work is of the form
$(-1)^{s+1} \partial_x^s b(x) \partial_x^s$ instead of
$(-1)^{s+1} \left( \partial_x a(x) \partial_x \right)^s$.

Multiplication of \eqref{eq:RHS-smart} with $\vec{1}^T \mat{M}$ results in
\begin{equation}
  - \epsilon \, \vec{1}^T \mat{D}[^T] \mat{M} \mat{a} \mat{D}
  \left( \mat{M}[^{-1}] \mat{D}[^T] \mat{M} \mat{a} \mat{D} \right)^{s-1} \vec{u}
  = 0,
\end{equation}
since the derivative is exact for constants. Therefore, the resulting scheme is
conservative across elements. 

Multiplying \eqref{eq:RHS-smart} by $\vec{u}^T \mat{M}$ and using the symmetry of
$\mat{M}$ yields
\begin{equation}
\begin{aligned}
  &
  - \epsilon \, \vec{u}^T \mat{D}[^T] \mat{M} \mat{a} \mat{D}
  \left( \mat{M}[^{-1}] \mat{D}[^T] \mat{M} \mat{a} \mat{D} \right)^{s-1} \vec{u}
  \\&
  =
  - \epsilon \, \vec{u}^T \mat{D}[^T] \mat{M} \mat{a} \mat{D} \mat{M}[^{-1}] \mat{M}
  \left( \mat{M}[^{-1}] \mat{D}[^T] \mat{M} \mat{a} \mat{D} \right)^{s-1} \vec{u}
  \\&
  = - \epsilon
    \left[
      \left( \mat{M}[^{-1}] \mat{D}[^T] \mat{a}[^T] \mat{M} \mat{D} \right) \vec{u}
    \right]^T
    \mat{M}
    \left[
      \left( \mat{M}[^{-1}] \mat{D}[^T] \mat{M} \mat{a} \mat{D} \right)^{s-1} \vec{u}
    \right].
\end{aligned}
\end{equation}
If the multiplication operator $\mat{a}$ is self-adjoint with respect to $\mat{M}$,
i.e. $\mat{M} \mat{a} = \mat{a}[^T] \mat{M}$, this becomes by induction
\begin{equation}
\begin{aligned}
  &
  - \epsilon \, \vec{u}^T \mat{D}[^T] \mat{M} \mat{a} \mat{D}
  \left( \mat{M}[^{-1}] \mat{D}[^T] \mat{M} \mat{a} \mat{D} \right)^{s-1} \vec{u}
  \\&
  = - \epsilon
    \left[
      \left( \mat{M}[^{-1}] \mat{D}[^T] \mat{M} \mat{a} \mat{D} \right) \vec{u}
    \right]^T
    \mat{M}
    \left[
      \left( \mat{M}[^{-1}] \mat{D}[^T] \mat{M} \mat{a} \mat{D} \right)^{s-1} \vec{u}
    \right]
  \\&
  =
  \begin{cases}
    - \epsilon
      \left[
        \left( \mat{M}[^{-1}] \mat{D}[^T] \mat{M} \mat{a} \mat{D} \right)^{s/2} \vec{u}
      \right]^T
      \mat{M}
      \left[
        \left( \mat{M}[^{-1}] \mat{D}[^T] \mat{M} \mat{a} \mat{D} \right)^{s/2} \vec{u}
      \right],
    & s \text{ even},
    \\
    - \epsilon
      \left[
        \left( \mat{M}[^{-1}] \mat{D}[^T] \mat{M} \mat{a} \mat{D} \right)^\frac{s-1}{2} \vec{u}
      \right]^T
      \mat{D}[^T] \mat{M} \mat{a} \mat{D}
      \left[
        \left( \mat{M}[^{-1}] \mat{D}[^T] \mat{M} \mat{a} \mat{D} \right)^\frac{s-1}{2} \vec{u}
      \right],
    & s \text{ odd}.
  \end{cases}
\end{aligned}
\end{equation}

If $a \geq 0$, this is non-positive for a nodal basis with diagonal mass matrix
$\mat{M}$, since the multiplication matrix
$\mat{a} = \diag{a(\xi_0), \dots, a(\xi_p)}$ is diagonal and
has non-negative entries. Thus, $\mat{a}$ is $\mat{M}$-self-adjoint and
$\mat{M} \mat{a} = \mat{M} \sqrt{\strut \mat{a}}^2 =
\sqrt{\strut \mat{a}} \mat{M} \sqrt{\strut \mat{a}}$,
since diagonal matrices commute. Thus, the resulting scheme is stable in the
discrete norm $\norm{\cdot}_M$ induced by $\mat{M}$.

If $\mat{M}$ represents the $L_2$ scalar product, $a \geq 0$ is a polynomial and
multiplication is given by an exact multiplication of polynomials followed by
an exact $L_2$ projection, $\mat{a}$ is $\mat{M}$-self-adjoint since for arbitrary
poylnomials $u, v$ of degree $\leq p$
\begin{equation}
  \vec{v}^T \mat{M} \mat{a} \vec{u}
  = \int v \, \proj(a u)
  = \int v  \, a  \, u
  = \int \proj(a v) \, u
  = \vec{v}^T \mat{a}[^T] \mat{M} \vec{u}.
\end{equation}
Additionally, $\mat{M} \mat{a}$ is positive semi-definite, since
\begin{equation}
  \vec{v}^T \mat{M} \mat{a} \vec{v}
  = \int v \, \proj(a v)
  = \int v  \, a  \, v
  = \int a \, v^2
  \geq 0
\end{equation}
for an arbitrary polynomial $v$ of degree $\leq p$. Therefore, the resulting scheme
is stable.
These results are summed up in the following
\begin{lem}
\label{lem:semidiscrete}
  Augmenting a conservative and stable SBP CPR method for the scalar conservation
  law \eqref{eq:scalar-CL}
  \begin{equation}
    \partial_t u + \partial_x f(u) = 0
  \end{equation}
  with the right hand side \eqref{eq:RHS-smart}
  \begin{equation}
    - \epsilon 
      \left( \mat{M}[^{-1}] \mat{D}[^T] \mat{M} \mat{a} \mat{D} \right)^s
      \vec{u},
  \end{equation}
  where $a|_{\Omega} \geq 0$ is a polynomial fulfilling $a|_{\partial \Omega} = 0$,
  results in a conservative and stable semidiscrete scheme if
  \begin{itemize}
    \item a nodal basis with diagonal norm matrix $\mat{M}$
    \item or a modal basis with exact $L_2$ norm and multiplication using
          exact $L_2$ projection
  \end{itemize}
  is used. Bases fulfilling this conditions are nodal bases using Gauß-Legendre
  or Lobatto-Legendre nodes (with lumped mass matrix) and a modal Legendre basis.
\end{lem}
Here, \emph{conservative} refers to conservation across elements and \emph{stable}
refers to stability in the discrete norm $\norm{\cdot}_M$ induced by $\mat{M}$,
approximating the $L_2$ norm.

\subsection{Eigenvalues of the discrete dissipation operator}

Choosing $a(x) = 1 - x^2$ for the standard element $\Omega = [-1,1]$,
the viscosity operator $\partial_x a(x) \partial_x$ for $s = 1$ on the right hand
side of \eqref{eq:scalar-CL-RHS} yields Legendre's differential equation
\begin{equation}
\label{eq:Legendre-eigenvalues}
  \partial_x \left( (1-x^2) \partial_x \phi_n(x) \right)
  = - n (n+1) \phi_n(x).
\end{equation}
Here, $n \in \N_0$ and $\phi_n$ is the Legendre polynomial of degree $n$. Thus,
the Legendre polynomials $\phi_n$ are eigenvectors of the continuous viscosity
operator with eigenvalues $- n (n+1)$.

In the discrete setting using polynomials of degree $\leq p$, this can only hold
for $n \leq p-1$, since $\phi_p$ has degree $p$ and $((1-x^2) \partial_x \phi_p(x)$
is represented by $\mat{a} \mat{D} \vec{\phi}_p$, i.e. a polynomial of degree
$\leq p$. Thus, $\mat{D} \mat{a} \mat{D} \vec{\phi}_p$ represents a polynomial
of degree $\leq p-1$ and especially not $\phi_p$.

Considering modal and nodal bases separately, the eigenvalues of the discrete
viscosity operator given by \eqref{eq:RHS-smart} for $s = 1$ (and thus for
arbitrary $s$ by multiplication) are computed in the following paragraphs.

\subsubsection*{Modal Legendre basis}

At first, a modal Legendre basis with exact $L_2$ scalar product is assumed.
Since multiplication with $a(x) = 1 - x^2$ increases the degree at most by 2,
the discretisation \eqref{eq:RHS-smart} yields the correct eigenvalues for
$n \in \set{0, \dots, p-1}$. Using the orthogonality of Legendre polynomials,
the $k$-th coefficient of any vector $\vec{v}$ is given by
$\norm{\vec{\phi}_k}_M^2 \left[ \vec{v} \right]_k = \vec{\phi}_k^T \mat{M} \vec{v}$.
Therefore,
\begin{equation}
  \norm{\vec{\phi}_k}_M^2  \left[ -\mat{M}[^{-1}] \mat{D}[^T] \mat{M} \mat{a} \mat{D} \vec{\phi}_n \right]_k
  = - \vec{\phi}_k^T \mat{D}[^T] \mat{M} \mat{a} \mat{D} \vec{\phi}_n.
\end{equation}
For $n \in \set{0, \dots, p-1}$, the right hand side is evaluated exactly since
no projection is necessary. Thus,
\begin{equation}
\begin{aligned}
  \norm{\vec{\phi}_k}_M^2 
  \left[ -\mat{M}[^{-1}] \mat{D}[^T] \mat{M} \mat{a} \mat{D} \vec{\phi}_n \right]_k
  &=
  - \int \partial_x \phi_k(x) \cdot a(x) \partial_x \phi_n(x)
  \\&=
  - n (n+1) \norm{\phi_n}_{L_2[-1,1]}^2 \delta_{k n},
\end{aligned}
\end{equation}
for $k \in \set{0, \dots, p}$. Here, $\delta_{k n} = 1$ for $k = n$ and
$\delta_{k n} = 0$ for $k \neq n$. This can simply be rewritten as
\begin{equation}
  -\mat{M}[^{-1}] \mat{D}[^T] \mat{M} \mat{a} \mat{D} \vec{\phi}_n
  = - n (n+1) \vec{\phi}_n, \quad n \in \set{0, \dots, p-1}.
\end{equation}

For $n = p$, equation \eqref{eq:RHS-Legendre-p} of the appendix can be used,
resulting in
\begin{equation}
  \proj \left( (1-x^2) \partial_x \phi_p(x) \right) = \frac{p (p+1)}{2p+1} \phi_{p-1}(x).
\end{equation}
Therefore,
\begin{equation}
  \mat{a} \mat{D} \vec{\phi}_p
  = \frac{p (p+1)}{2p+1} \vec{\phi}_{p-1}.
\end{equation}
Using equation \eqref{eq:Legendre-derivative} of the appendix, the derivative
of a Legendre polynomial is given by
\begin{equation}
  \mat{D} \vec{\phi}_k
  = (2k-1) \vec{\phi}_{k-1} + \mat{D} \vec{\phi}_{k-2}
  = (2k-1) \vec{\phi}_{k-2} + (2k-5) \vec{\phi}_{k-3} + \dots .
\end{equation}
Therefore, by using the orthogonality of Legendre polynomials, for
$k \in \set{0, \dots, p-1}$,
\begin{equation}
  \norm{\vec{\phi}_k}_M^2  \left[
    -\mat{M}[^{-1}] \mat{D}[^T] \mat{M} \mat{a} \mat{D} \vec{\phi}_p \right]_k
  = - \vec{\phi}_k^T \mat{D}[^T] \mat{M} \mat{a} \mat{D} \vec{\phi}_p = 0.
\end{equation}
Additionally, for $k = p$,
\begin{equation}
\begin{aligned}
  &
  \norm{\vec{\phi}_p}_M^2  \left[
  -\mat{M}[^{-1}] \mat{D}[^T] \mat{M} \mat{a} \mat{D} \vec{\phi}_p \right]_p
  =
  - \vec{\phi}_p^T \mat{D}[^T] \mat{M} \mat{a} \mat{D} \vec{\phi}_p
  \\=&
  - (2p-1) \frac{p (p+1)}{2p+1} \vec{\phi}_{p-1}^T \mat{M} \vec{\phi}_{p-1}
  =
  - (2p-1) \frac{p (p+1)}{2p+1} \norm{\vec{\phi}_{p-1}}_M^2
  \\=&
  - (2p-1) \frac{p (p+1)}{2p+1} \norm{\phi_{p-1}}_{L_2[-1,1]}^2
  =
  - (2p-1) \frac{p (p+1)}{2p+1} \frac{2}{2p-1}
  \\=&
  - p (p+1) \frac{2}{2p+1}
  =
  - p (p+1) \norm{\vec{\phi}_p}_M^2.
\end{aligned}
\end{equation}
Thus, for a modal Legendre basis, the viscosity operator \eqref{eq:RHS-smart}
yields the correct eigenvalues for the Legendre polynomials
\begin{equation}
  -\mat{M}[^{-1}] \mat{D}[^T] \mat{M} \mat{a} \mat{D} \vec{\phi}_n
  = - n (n+1) \vec{\phi}_n, \quad n \in \set{0, \dots, p}.
\end{equation}

\subsubsection*{Nodal bases}

Considering a nodal basis with diagonal norm matrix $\mat{M}$, corresponding to
some quadrature rule with positive weights, the order of the quadrature is
important. If the quadrature given by $\mat{M}$ is exact for polynomials of
degree $\leq q$,
\begin{equation}
  \norm{\vec{\phi}_k}_M^2  \left[
    -\mat{M}[^{-1}] \mat{D}[^T] \mat{M} \mat{a} \mat{D} \vec{\phi}_n \right]_k
  = - \vec{\phi}_k^T \mat{D}[^T] \mat{M} \mat{a} \mat{D} \vec{\phi}_n
  , \quad 0 \leq k + p \leq q,
\end{equation}
since $\partial_x \phi_k \, a \, \partial_x \phi_n$ is of degree $\leq q$.

For Gauß-Legendre nodes, corresponding to a quadrature of degree $2p+1$, this yields
\begin{equation}
  -\mat{M}[^{-1}] \mat{D}[^T] \mat{M} \mat{a} \mat{D} \vec{\phi}_n
  = - n (n+1) \vec{\phi}_n, \quad n \in \set{0, \dots, p}.
\end{equation}
Lobatto-Legendre nodes result in a quadrature of degree $2p-1$, thus 
\begin{equation}
  -\mat{M}[^{-1}] \mat{D}[^T] \mat{M} \mat{a} \mat{D} \vec{\phi}_n
  = - n (n+1) \vec{\phi}_n, \quad n \in \set{0, \dots, p-1}.
\end{equation}
By the same reason, for $n=p$, 
\begin{equation}
  \norm{\vec{\phi}_k}_M^2  \left[
    -\mat{M}[^{-1}] \mat{D}[^T] \mat{M} \mat{a} \mat{D} \vec{\phi}_p \right]_k
  = - \vec{\phi}_k^T \mat{D}[^T] \mat{M} \mat{a} \mat{D} \vec{\phi}_p
  = 0, \quad k \in \set{0, \dots, p-1}.
\end{equation}
The case $k = n = p$ is a bit more complicated. Using equation
\eqref{eq:RHS-Legendre-p} of the appendix,
\begin{equation}
  (1-x^2) \partial_x \phi_p(x)
  = \frac{p (p+1)}{2p+1} \left( \phi_{p-1}(x) - \phi_{p+1}(x) \right).
\end{equation}
Again, by equation \eqref{eq:Legendre-derivative},
\begin{equation}
  \partial_x \phi_p
  = (2p-1) \phi_{p-1} + (2p-5) \phi_{p-3} + \dots .
\end{equation}
Since Lobatto-Legendre quadrature is exact for polynomials of degree $\leq 2p-1$
and Legendre polynomials are orthogonal,
\begin{equation}
  - \vec{\phi}_p^T \mat{D}[^T] \mat{M} \mat{a} \mat{D} \vec{\phi}_p
  = - (2p-1) \frac{p (p+1)}{2p+1} \vec{\phi}_{p-1}^T \mat{M}
    \left( \vec{\phi}_{p-1} - \vec{\phi}_{p+1} \right).
\end{equation}
Inserting the equations \eqref{eq:norm-lobatto-phi_p} and
\eqref{eq:lobatto-phi_p-1-phi_p+1} of the appendix finally results in
\begin{equation}
  - \vec{\phi}_p^T \mat{D}[^T] \mat{M} \mat{a} \mat{D} \vec{\phi}_p
  = - p (p+1) \frac{2p-1}{2p+1} \left(
    \frac{2}{p} - \frac{2}{p} \frac{ p }{ 2p-1 }
    \right).
\end{equation}
The term in brackets can be simplified as
\begin{equation}
  \frac{2}{p} \left( 1 - \frac{p}{2p-1} \right)
  = \frac{ 2 (p-1) }{ p (2p-1) }
  = \frac{2}{2p-1} \frac{p-1}{p}.
\end{equation}
This term is zero for $p = 1$ and positive for $p > 1$. Additionally, it is
monotonically decreasing for $p > 1$ and bounded from above by $2 / (2p-1)$.

\subsubsection*{Summary}

These results are summed up in the following
\begin{lem}
\label{lem:eigenvalues}
  The discrete viscosity operator $- \mat{M}[^{-1}] \mat{D}[^T] \mat{M} \mat{a} \mat{D} u$
  of the right hand side \eqref{eq:RHS-smart}
  \begin{itemize}
    \item 
    has the same eigenvalues $- n (n+1)$ for the Legendre polynomials
    $\phi_n, n \in \set{0, \dots, p}$ as the continuous operator
    $\partial_x (1-x^2) \partial_x$, if a modal Legendre or a nodal Gauß-Legendre
    basis is used.
    
    \item
    has the same eigenvalues $- n (n+1)$ for the Legendre polynomials
    $\phi_n, n \in \set{0, \dots, p-1}$ as the continuous operator
    $\partial_x (1-x^2) \partial_x$, if a nodal Lobatto-Legendre basis is
    used. The eigenvalue for $\phi_p$ is non-positive and bigger than $- p (p+1)$,
    i.e. the viscosity operator yields less dissipation of the highest
    mode compared to the exact value, obtained by modal Legendre and nodal
    Gauß-Legendre bases.
  \end{itemize}
\end{lem}

\subsection{Discrete setting}

In order to get a working numerical scheme, a time discretisation has to be
introduced. For simplicity, an explicit Euler method will be considered. Thus,
the development in the standard element during one time step $\Delta t$ is
given by
\begin{equation}
\label{eq:exp-Euler}
  \vec{u}
  \mapsto \vec{u}_+ := \vec{u} + \Delta t \, \partial_t \vec{u}.
\end{equation}
However, if the fully discrete scheme using an explicit Euler method is both
conservative and stable, then a \emph{strong-stability preserving} (SSP)
has the same properties, since it consists of a convex combination of Euler
steps, see inter alia the monograph \citex{}{by} \citet{gottlieb2011strong} and references
cited therein.

Using an SBP CPR semidiscretisation to compute the time derivative
$\partial_t \vec{u}$ for a scalar conservation law \eqref{eq:scalar-CL} without
artificial viscosity term, the norm after one Euler step is given by
\begin{equation}
\begin{aligned}
  \norm{ \vec{u}_+ }_M^2
  = \vec{u}_+^T \mat{M} \vec{u}_+
  &
  = \vec{u}^T \mat{M} \vec{u}
    + 2 \Delta t \, \vec{u}^T \mat{M} \partial_t \vec{u}
    + (\Delta t)^2 \partial_t \vec{u}^T \mat{M} \partial_t \vec{u}
  \\&
  = \norm{ \vec{u} }_M^2
    + 2 \Delta t \, \scp{ \vec{u} }{ \partial_t \vec{u} }_M
    + (\Delta t)^2 \norm{ \partial_t \vec{u} }_M^2.
\end{aligned}
\end{equation}
Here, the second term on the right hand side
$2 \Delta t \, \scp{ \vec{u} }{ \partial_t \vec{u} }_M$
has been estimated for the semidiscretisation, yielding only boundary terms that
can be controlled by the numerical flux and results consequently in a stable
scheme.

However, the last term $(\Delta t)^2 \norm{ \partial_t \vec{u} }_M^2 \geq 0$
is non-negative and adds an undesired effect, i.e. increases the norm and may
trigger instabilities. Thus, the basic idea is to add artificial dissipation and
choose the parameters appropriately in order to damp the undesired energy growth.

Assuming a fixed function $a$ and order $s$, the strength $\epsilon$ can be
estimated in the following way. Denoting the time derivative obtained by the
underlying SBP CPR method without artificial dissipation by $\partial_t \vec{u}$,
introducing the artificial viscosity term \eqref{eq:RHS-smart} with strength
$\epsilon$ yields
\begin{equation}
\label{eq:av-derivative}
  \partial_t \vec{u}^\epsilon
  = 
  \partial_t \vec{u}
  - \epsilon \left( \mat{M}[^{-1}] \mat{D}[^T] \mat{M} \mat{a} \mat{D} \right)^s \vec{u}.
\end{equation}
Abbreviating the artificial dissipation operator \eqref{eq:RHS-smart} as
\begin{equation}
  \mat{A}[^s]
  :=
  \left(\mat{M}[^{-1}] \mat{D}[^T] \mat{M} \mat{a} \mat{D} \right)^s,
\end{equation}
the norm after one explicit Euler step with artificial dissipation is
\begin{equation}
\label{eq:av-euler-condition}
\begin{aligned}
  \norm{ \vec{u}^\epsilon_+ }_M^2
  =&
  \norm{ \vec{u} }_M^2
  + 2 \Delta t \scp{ \vec{u} }{ \partial_t \vec{u}^\epsilon }_M
  + (\Delta t)^2 \norm{ \partial_t \vec{u}^\epsilon }_M^2
  \\=&
  \norm{ \vec{u} }_M^2
  + 2 \Delta t \scp{ \vec{u} }{ \partial_t \vec{u} }_M
  - 2 \epsilon \Delta t \scp{\vec{u}}{\mat{A}[^s] \vec{u}}_M
  + (\Delta t)^2 \norm{ \partial_t \vec{u}^\epsilon }_M^2
\end{aligned}
\end{equation}
Again, $\scp{ \vec{u} }{ \partial_t \vec{u} }_M$ can be estimated in terms of
boundary values and numerical fluxes. Thus, the two last terms shall cancel out,
resulting in an estimate similar to the semidiscrete one, i.e.
\begin{equation}
  \norm{ \vec{u}^\epsilon_+ }_M^2
  =
  \norm{ \vec{u} }_M^2 
  + 2 \Delta t \scp{ \vec{u} }{ \partial_t \vec{u} }_M.
\end{equation}
Then, the fully discrete scheme will be both conservative across elements and
stable.

Using \eqref{eq:av-euler-condition}, this condition can be rewritten as
\begin{equation}
\begin{aligned}
 0 
 =&
 - 2 \epsilon \scp{\vec{u}}{\mat{A}[^s] \vec{u}}_M 
 + \Delta t \norm{ \partial_t \vec{u}^\epsilon }_M^2
 \\
 =&
 - 2 \epsilon \scp{\vec{u}}{\mat{A}[^s] \vec{u}}_M 
 + \Delta t \left( \norm{ \partial_t \vec{u} }_M^2
 - 2 \epsilon \scp{\partial_t \vec{u}}{\mat{A}[^s] \vec{u}}_M 
 + \epsilon^2 \norm{ \mat{A}[^s] \vec{u} }_M^2 \right),
\end{aligned}
\end{equation}
which is equivalent to
\begin{equation}
\label{eq:quadratic}
  \epsilon^2 \underbrace{\left( \Delta t \norm{ \mat{A}[^s] \vec{u} }_M^2 \right)}_{=: A} 
  + \epsilon \underbrace{\left( - 2 \scp{\vec{u}}{\mat{A}[^s] \vec{u}}_M 
    - 2 \Delta t \scp{\partial_t \vec{u}}{\mat{A}[^s] \vec{u}}_M \right)}_{=: B}
  + \underbrace{\left( \Delta t \norm{ \partial_t \vec{u} }_M^2 \right)}_{=: C}
  = 0.
\end{equation}
The (possibly complex) roots of this equation for $A \neq 0$ are given by
\begin{equation}
  \epsilon_{1/2} = \frac{1}{2 A} \left( -B \pm \sqrt{ B^2 -4 A C } \right).
\end{equation}
Since for a sufficiently small time step $\Delta t$ the discriminant $B^2 - 4 A C$ 
is non negative if the solution is not constant, there is at least one real
solution $\epsilon$. Additionally, both $-B$ and $A C$ are positive for
sufficiently small $\Delta t$, since the artificial dissipation operator
$\mat{A}$ is positive semi-definite, i.e.
\begin{equation}
\label{eq:av-condition-Delta-t}
  B^2 - 4 A C > 0
  ,\quad
  -B          > 0
  ,\qquad \text{ if } \Delta t \text{ is small enough and } \mat{A}[^s] \vec{u} \neq 0.
\end{equation}
Thus,
\begin{equation}
\begin{aligned}
  \epsilon_1
  \geq &
  \epsilon_2 
  = \frac{1}{2 A} \left( -B - \sqrt{ B^2 -4 A C } \right)
  \\
  \geq & \frac{1}{2 A} \left( -B + \sqrt{ B^2 } \right)  
  = 0,
\end{aligned}
\end{equation}
and the roots of the quadratic equation \eqref{eq:quadratic} are non-negative.
These results are summed up in the following
\begin{lem}
\label{lem:discrete}
  If a conservative and stable SBP CPR method for a scalar conservation law
  \eqref{eq:scalar-CL}
  \begin{equation}
    \partial_t u + \partial_x f(u) = 0
  \end{equation}
  is augmented with the artificial dissipation \eqref{eq:RHS-smart}
  \begin{equation}
    - \epsilon \left( \mat{M}[^{-1}] \mat{D}[^T] \mat{M} \mat{a} \mat{D} \right)^s \vec{u}
  \end{equation}
  on the right hand side, the fully discrete scheme using an explicit Euler
  method as time discretisation is both conservative and stable if
  \begin{itemize}
    \item 
    a nodal Gauß-Legendre / Lobatto-Legendre or a modal Legendre basis is used,
    
    \item
    $\scp{ \vec{u} }{ \mat{A}[^s] \vec{u} } > 0$, which will be fulfilled for
    the choice of $a$ described below if the solution $\vec{u}$ is not constant,
    
    \item
    the time step $\Delta t$ is small enough such that \eqref{eq:av-condition-Delta-t}
    is fulfilled,
    
    \item
    and the strength $\epsilon > 0$ is big enough.
  \end{itemize}
  
  If the other conditions are fulfilled, $\epsilon$ has to obey
  \begin{equation}
  \label{eq:av-epsilon-adaptive}
    \epsilon
    \geq \epsilon_2
    = \frac{1}{2 A} \left( -B - \sqrt{ B^2 -4 A C } \right),
  \end{equation}
  where $A, B,$ and $C$ from equation \eqref{eq:quadratic} are used.
\end{lem}

In numerical computations, the second (smaller) root $\epsilon_2$ is used as
strength and results in methods with highly desired stability properties,
as described in the next section.

However, it remains an interesting and yet unanswered question how to interpret
the existence of an additional solution $\epsilon_1$. Since this solution yields
a bigger strength, the resulting methods show higher dissipation, which might
be undesired in elements without discontinuities or for long time simulations.

Additionally, equation \eqref{eq:av-condition-Delta-t} limits the maximal time
step. This could be used for an adaptive control of the step size. This adaptive
strategy will be a topic of further research, together with more sophisticated
integration schemes. Thus, a simple limiting strategy is used for the numerical
experiments, i.e. if the time step is not small enough and equation
\eqref{eq:av-condition-Delta-t} is not fulfilled, the strength $\epsilon$
computed from \eqref{eq:av-epsilon-adaptive} might be negative. In this case,
to avoid instabilities, $\epsilon$ is set to zero, i.e. no artificial viscosity
is used in the corresponding elements. This phenomenon is strongly connected
with stability requirements of the viscous operator. Considering an explicit
Euler step for the equation $\partial_t \vec{u} = - \epsilon \mat{A}[^s] \vec{u}$,
the norm after one time step obeys
\begin{equation}
  \norm{\vec{u}_+}_M^2
  =
  \norm{\vec{u}}_M^2
  - 2 \, \epsilon \, \Delta t \scp{ \vec{u} }{ \mat{A}[^s] \vec{u} }_M
  + \epsilon^2 (\Delta t)^2 \norm{\mat{A}[^s] \vec{u}}_M^2.
\end{equation}
Therefore, in order to guarantee $\norm{\vec{u}_+}_M^2 \leq \norm{\vec{u}}_M^2$,
for $\mat{A}[^s] \vec{u} \neq 0$, $\Delta t$ is limited by
\begin{equation}
  \Delta t
  \leq
  \frac{ 2 \scp{ \vec{u} }{ \mat{A}[^s] \vec{u} }_M }{ \epsilon \norm{\mat{A}[^s] \vec{u}}_M^2 }.
\end{equation}

Since Lemma \ref{lem:semidiscrete} requires $a \big|_{[-1,1]} \geq 0$ to be a
polynomial fulfilling $a(\pm 1) = 0$, a simple choice is $a(x) = 1 - x^2$.
By this choice, the continuous artificial dissipation operators is related
to the eigenvalue equation of Legendre polynomials as described in the previous
section. Resulting implications
and connections with modal filtering are investigated \citex{in}{by}
\citet{glaubitz2016enhancing} in the second part of this article.

  \section{Numerical results}
\label{sec:numerical-results}

In order to augment the theoretical considerations of the previous chapters,
numerical experiments with and without artificial dissipation are presented in
this section.

\subsection{Linear advection with smooth initial condition}

Here, a numerical solution of the linear advection equation with constant velocity
\eqref{eq:lin-adv}
\begin{equation}
  \partial_t u + \partial_x u = 0
  ,\qquad
  u(0, x) = u_0(x) = \exp \left( -20 (x-1)^2 \right)
\end{equation}
with $N = 8$ elements using a Gauß-Legendre nodal basis of degree $\leq p = 7$
is computed in the domain $[0,2]$, equipped with periodic boundary conditions.
The time integration is performed by an explicit Euler method using $12 \cdot 10^4$
steps in the time interval $[0, 10]$ and a central numerical flux
$\fnum(u_-,u_+) = (u_- + u_+) / 2$ has been chosen for the semidiscretisation
\eqref{eq:lin-adv-semidisc}.

\begin{figure}[!htb]
  \centering
  \begin{subfigure}[b]{0.49\textwidth}
    \includegraphics[width=\textwidth]{%
      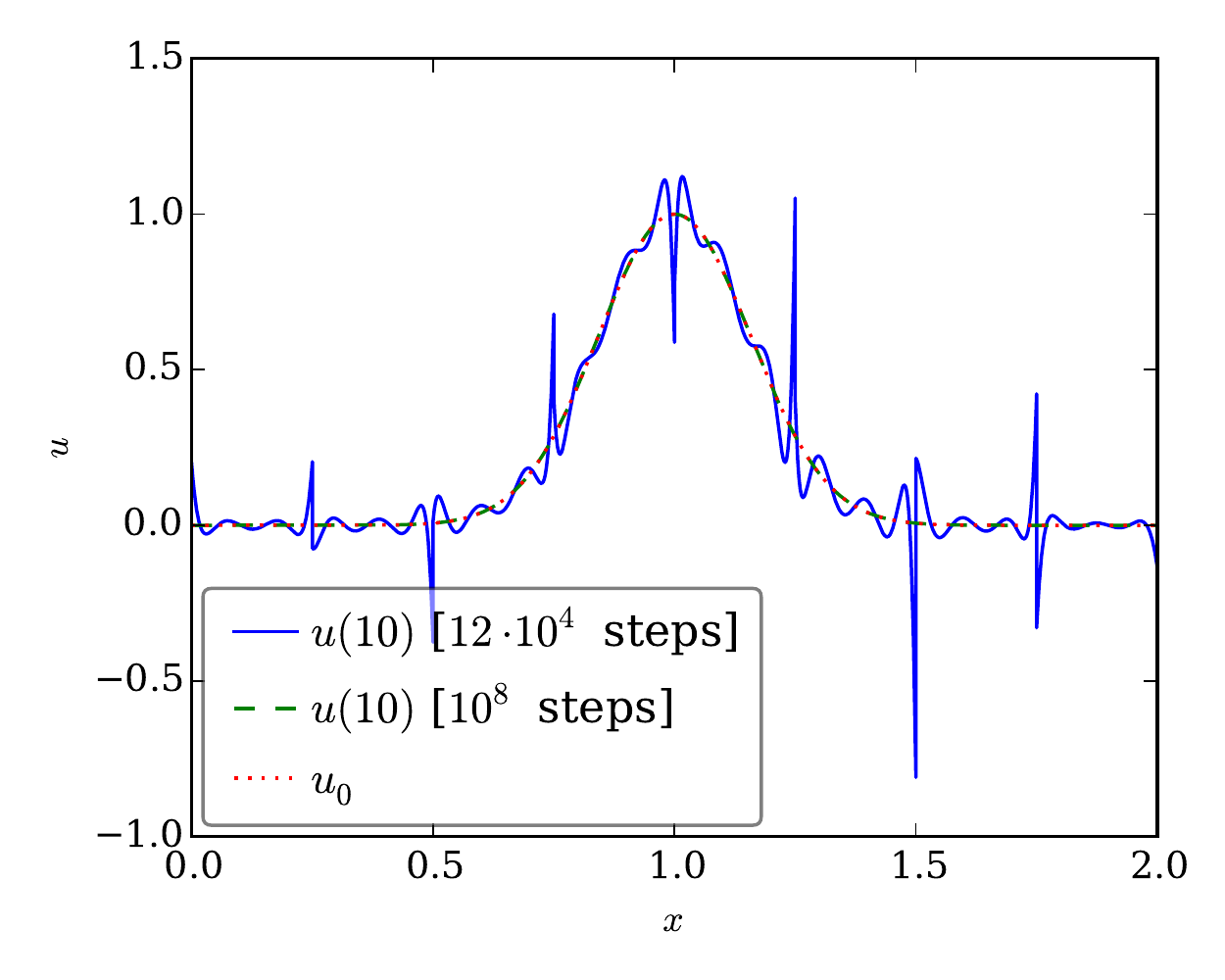}
    \caption{Solution computed without additional artificial dissipation.}
    \label{fig:linear_gauss_no_filtering_solution}
  \end{subfigure}%
  ~
  \begin{subfigure}[b]{0.49\textwidth}
    \includegraphics[width=\textwidth]{%
      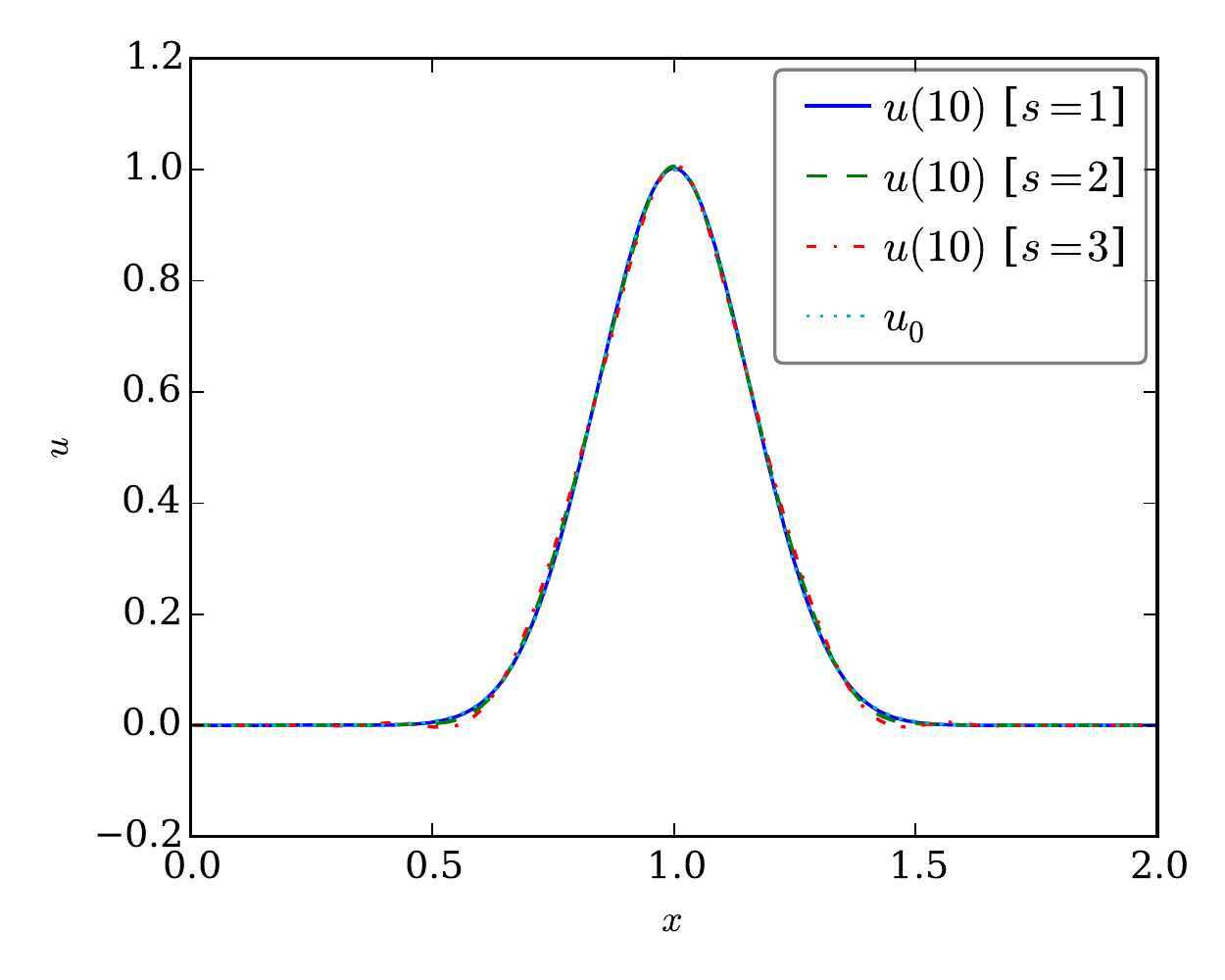}
    \caption{Solution computed using adaptive artificial dissipation.}
    \label{fig:linear_gauss_ad_adaptive_solution}
  \end{subfigure}%
  \\
  \begin{subfigure}[b]{0.49\textwidth}
    \includegraphics[width=\textwidth]{%
      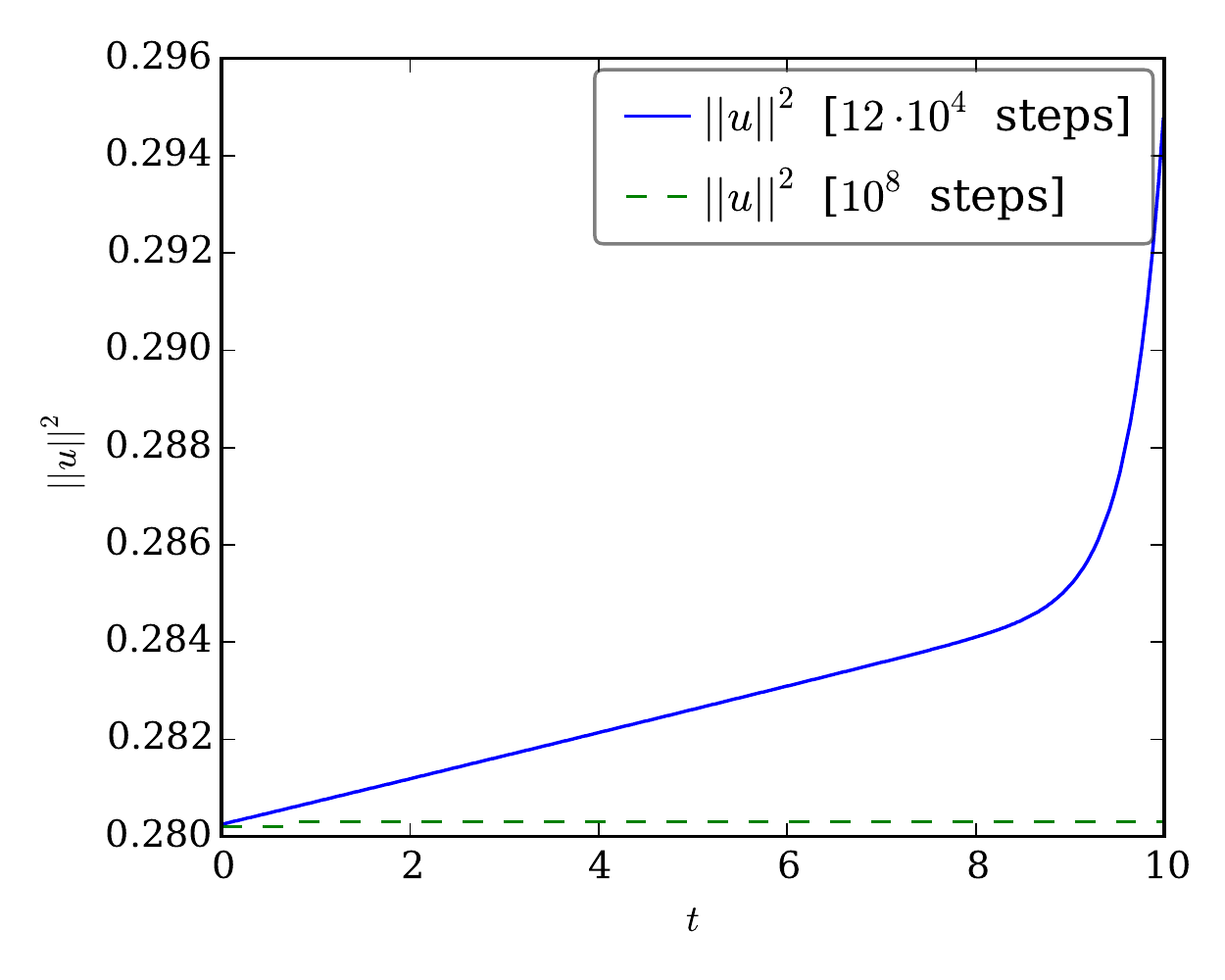}
    \caption{Energy of the solution computed without additional artificial dissipation.}
    \label{fig:linear_gauss_no_filtering_energy}
  \end{subfigure}%
  ~
  \begin{subfigure}[b]{0.49\textwidth}
    \includegraphics[width=\textwidth]{%
      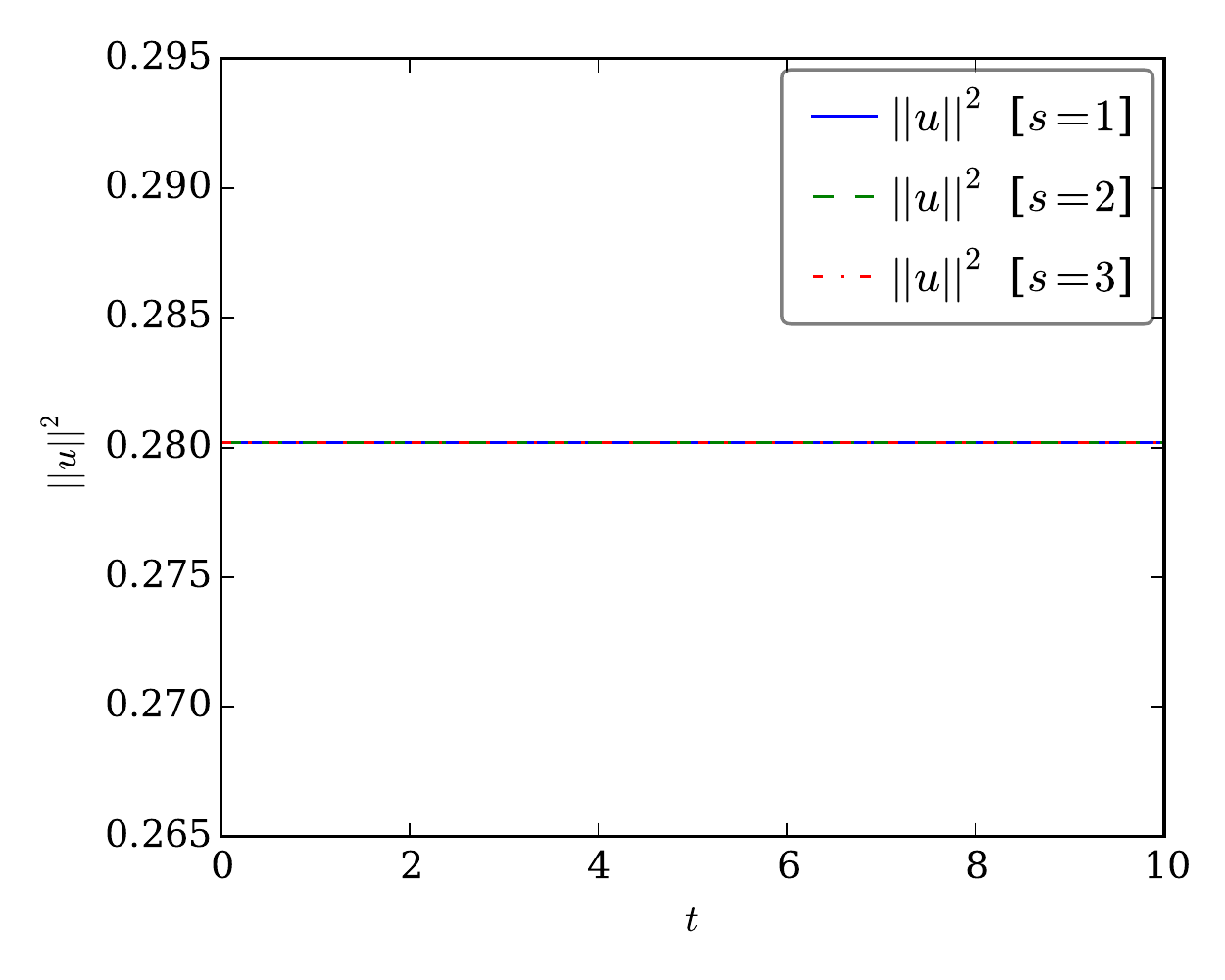}
    \caption{Energy of the solution computed using adaptive artificial dissipation.}
    \label{fig:linear_gauss_ad_adaptive_energy}
  \end{subfigure}%
  \caption{Numerical results for linear advection using $N = 8$ elements with
           polynomials of degree $\leq p = 7$.
           On the left hand side, no additional artificial dissipation has been
           used, whereas adaptive spectral viscosity of orders $s=1$, $s=2$, and
           $s=3$ has been used for the right hand side.
           }
  \label{fig:linear_gauss}
\end{figure}

As can be seen in Figure \ref{fig:linear_gauss_no_filtering_energy}, the energy
of the solution using $12 \cdot 10^4$ time steps is increasing, as expected.
This yields undesired oscillations in Figure \ref{fig:linear_gauss_no_filtering_solution}.
However, increasing the number of time steps to $10^8$ reduces the additional
term of order $(\Delta t)^2$. Therefore, the energy does not increase that much
and the solution has the desired smooth form.

However, the same effect can be achieved by adaptive artificial dissipation
using the estimate for the strength $\epsilon$ of Lemma \ref{lem:discrete}.
Using orders $s \in \set{1, 2, 3}$, the energy in Figure
\ref{fig:linear_gauss_ad_adaptive_energy} remains constant and the solutions in
Figure \ref{fig:linear_gauss_ad_adaptive_solution} look as expected. However,
there is a slight perturbation for $s = 3$ around $x \approx 1 \pm 0.5$.

\begin{figure}[!htb]
  \centering
  \begin{subfigure}[b]{0.49\textwidth}
    \includegraphics[width=\textwidth]{%
      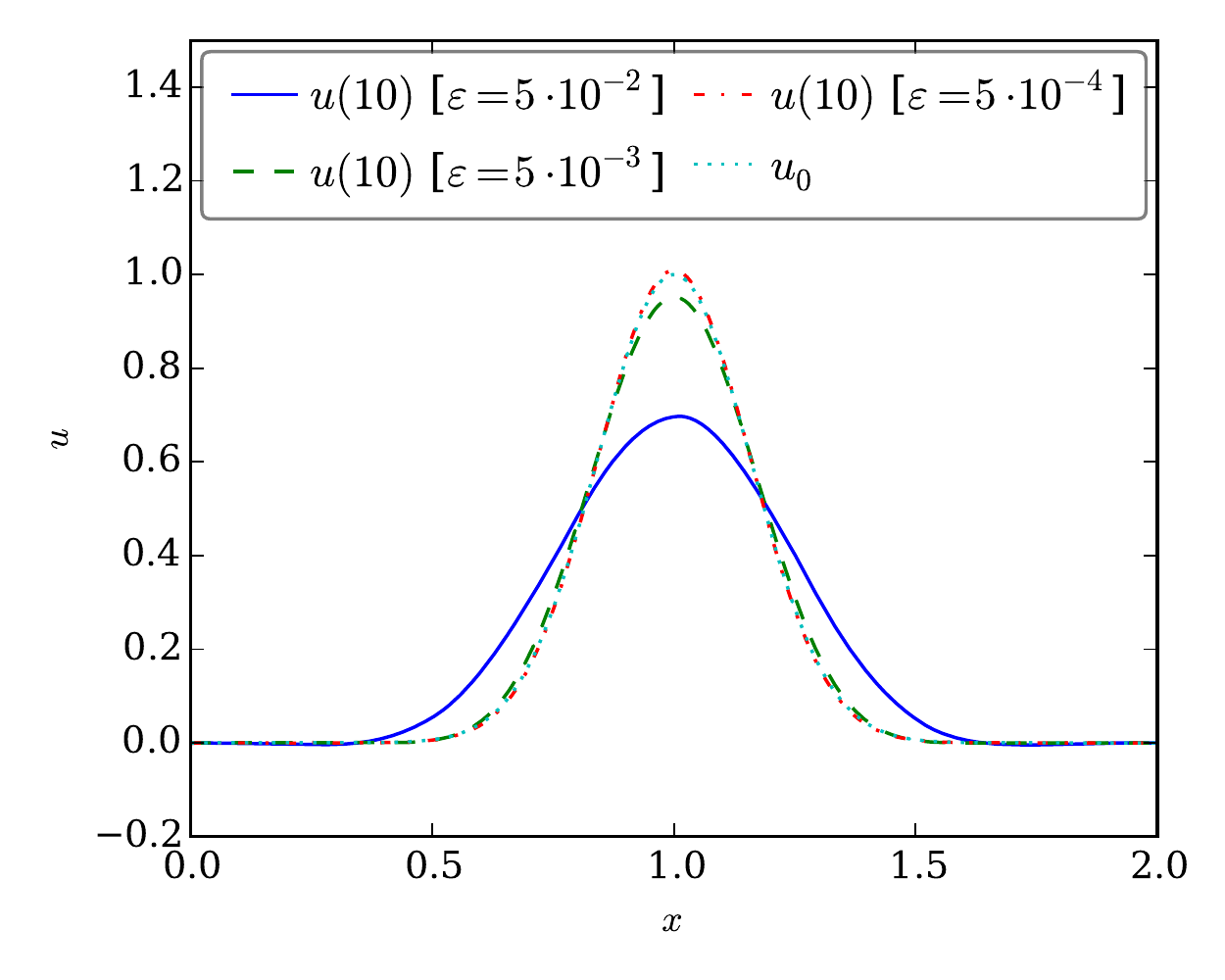}
    \caption{Solutions for order $s = 1$.}
  \end{subfigure}%
  ~
  \begin{subfigure}[b]{0.49\textwidth}
    \includegraphics[width=\textwidth]{%
      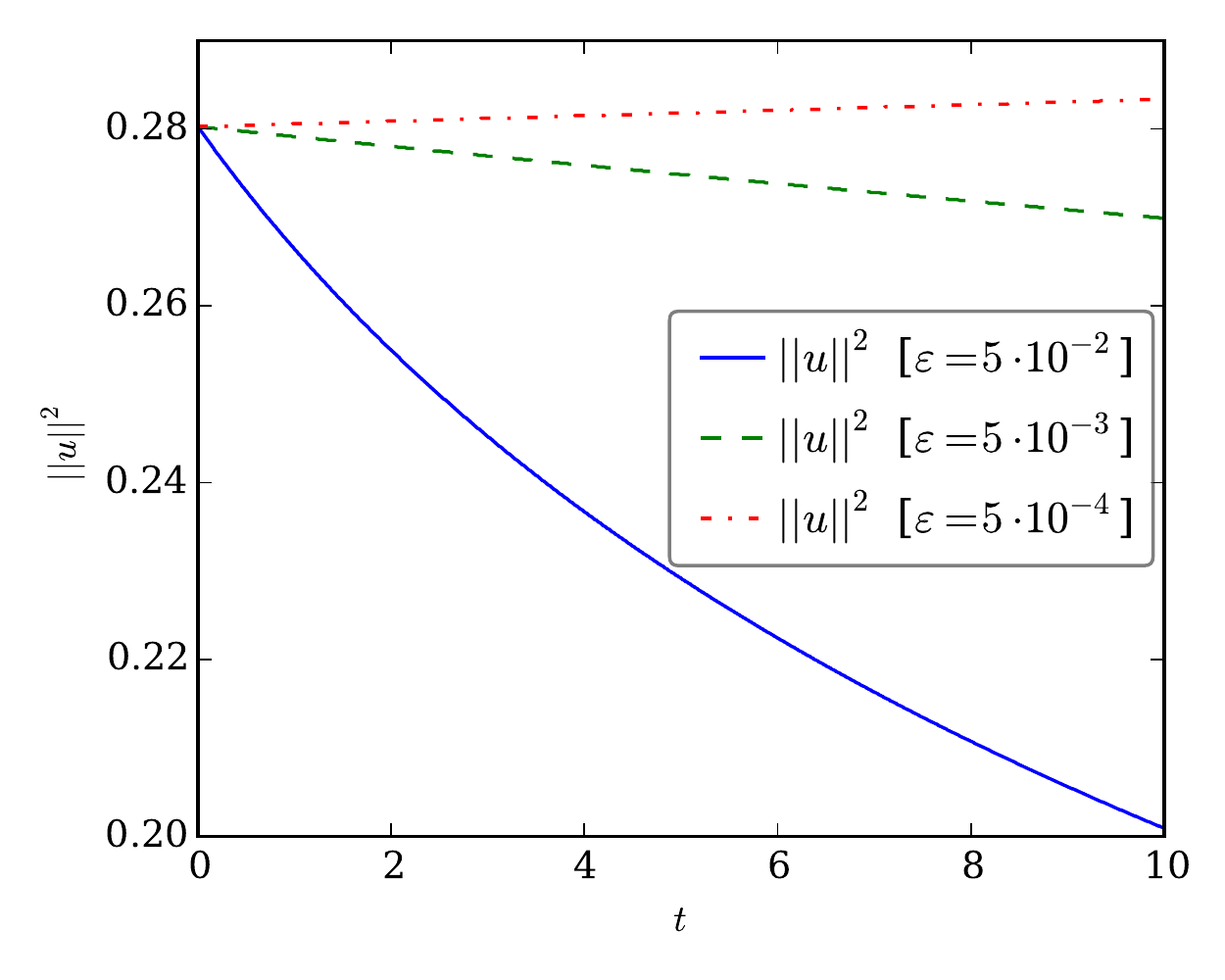}
    \caption{Energies for order $s = 1$.}
  \end{subfigure}%
  \\
  \begin{subfigure}[b]{0.49\textwidth}
    \includegraphics[width=\textwidth]{%
      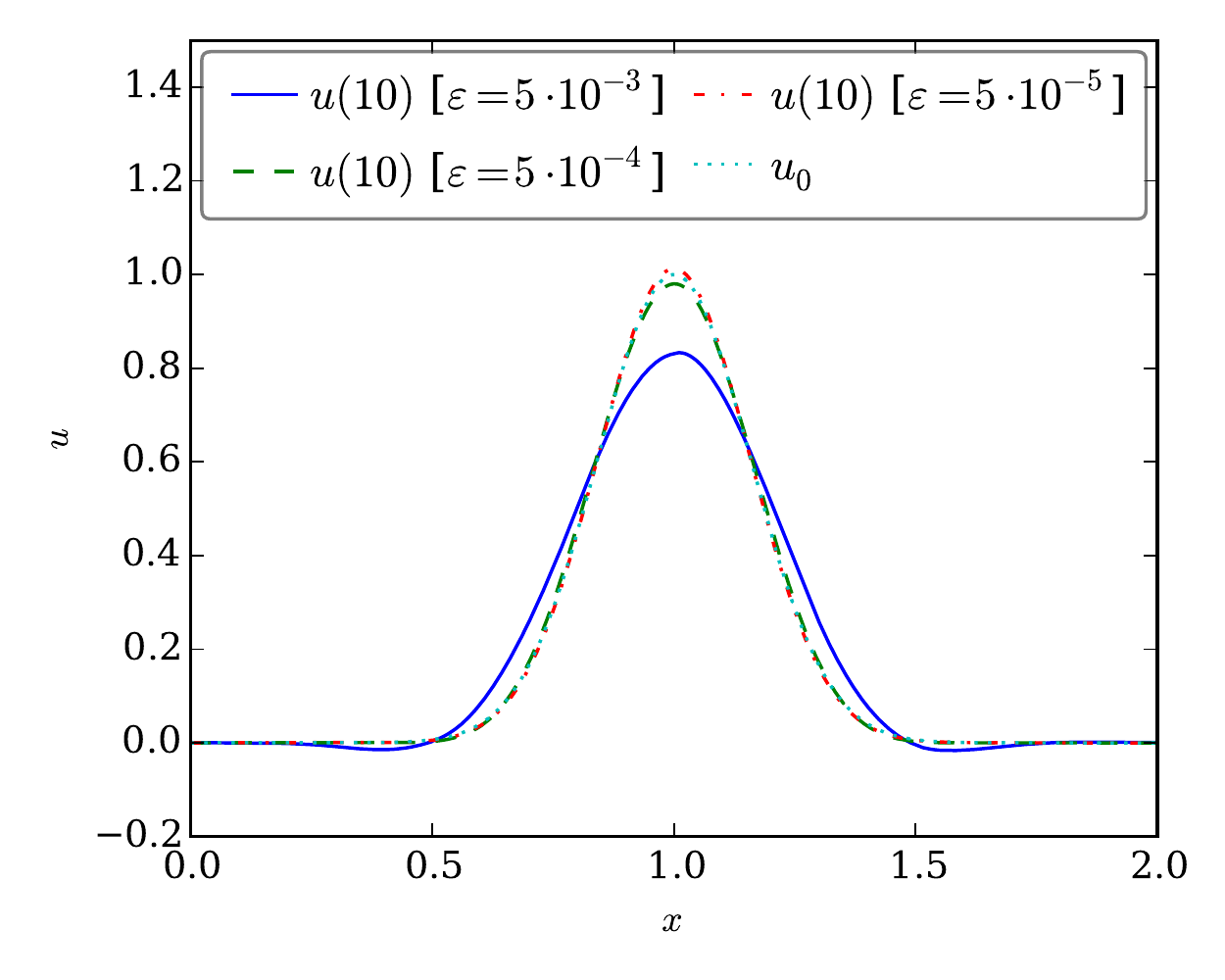}
    \caption{Solutions for order $s = 2$.}
  \end{subfigure}%
  ~
  \begin{subfigure}[b]{0.49\textwidth}
    \includegraphics[width=\textwidth]{%
      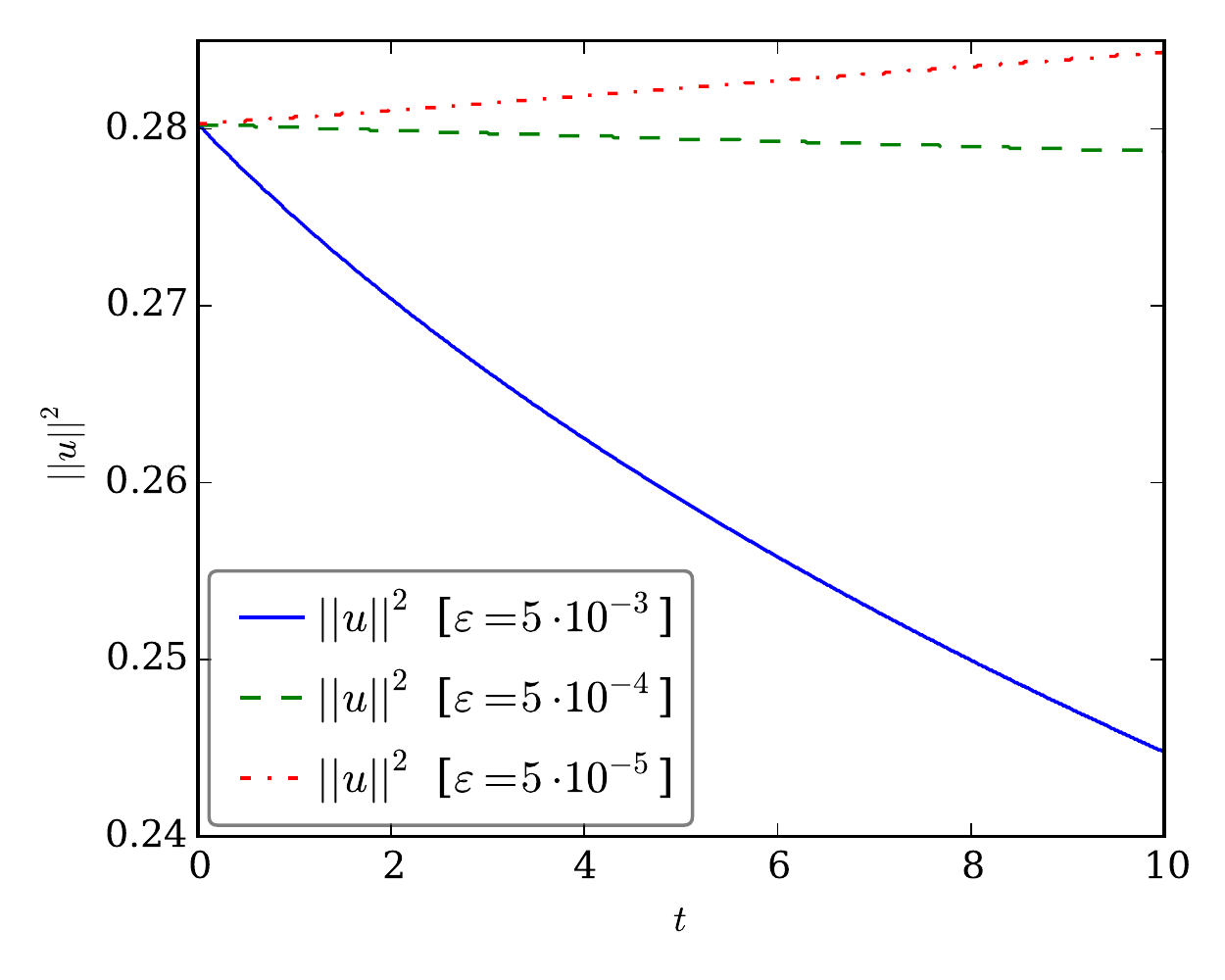}
    \caption{Energies for order $s = 2$.}
  \end{subfigure}%
  \caption{Numerical results for linear advection using $N = 8$ elements with
           polynomials of degree $\leq p = 7$ and constant artificial dissipation
           of orders $s \in \set{1,2}$ with various strengths $\epsilon$.
           On the left hand side, the solutions $u$ are shown, accompanied by the
           corresponding energies $\norm{u}^2$ on the right hand side.
           }
  \label{fig:linear_gauss_comparison}
\end{figure}

As can be seen in Figure \ref{fig:linear_gauss_comparison}, simple artificial
dissipation of fixed strength $\epsilon$ has a stabilising effect. The dissipation
of energy $\norm{u}^2$ is increasing with increasing order $s$ and strength
$\epsilon$, respectively. However, to get an acceptable result requires lengthy
experiments and fine tuning of the parameters by hand. Therefore, the adaptive
strategy of Lemma \ref{lem:discrete} provides an excellent alternative.

\subsection{Linear advection with discontinuous initial condition}

In order to see the influence in the presence of discontinuities, the linear
advection equation \eqref{eq:lin-adv}
\begin{equation}
  \partial_t u + \partial_x u = 0
  ,\qquad
  u(0, x) = u_0(x) =
  \begin{cases}
    1, & x \in [0.5, 1],\\
    0, & \text{otherwise},
  \end{cases}
\end{equation}
with periodic boundaries in the domain $[0,2]$ has been investigated during the
time interval $[0, 8]$. The semidiscretisation \eqref{eq:lin-adv-semidisc}
is used with an upwind numerical flux $\fnum(u_-,u_+) = u_-$ and rendered fully
discrete by an explicit Euler method.

\begin{figure}[!ht]
  \centering
  \begin{subfigure}[b]{0.49\textwidth}
    \includegraphics[width=\textwidth]{%
      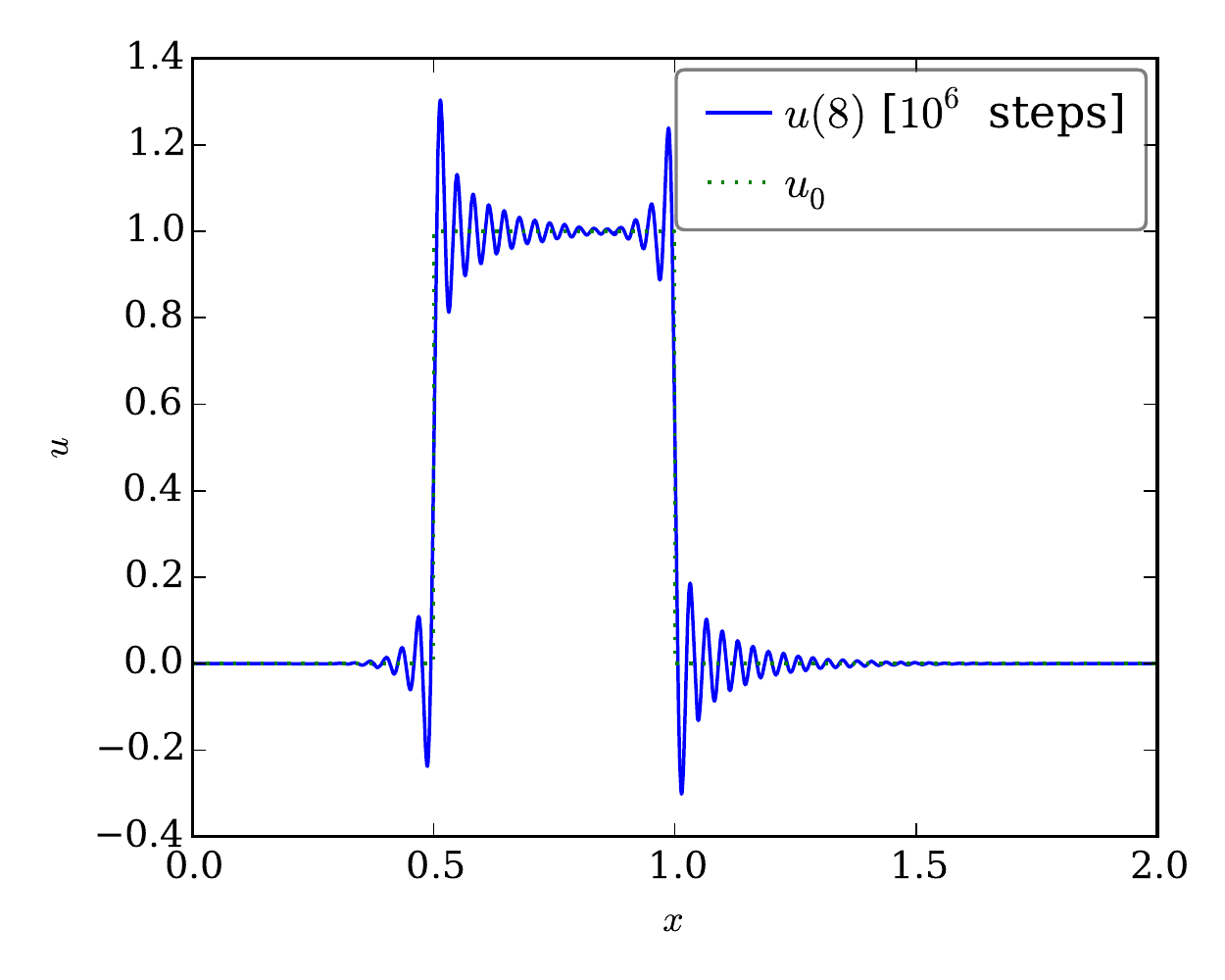}
    \caption{Solution computed without additional artificial dissipation.}
    \label{fig:linear_jump_no_filtering_solution}
  \end{subfigure}%
  ~
  \begin{subfigure}[b]{0.49\textwidth}
    \includegraphics[width=\textwidth]{%
      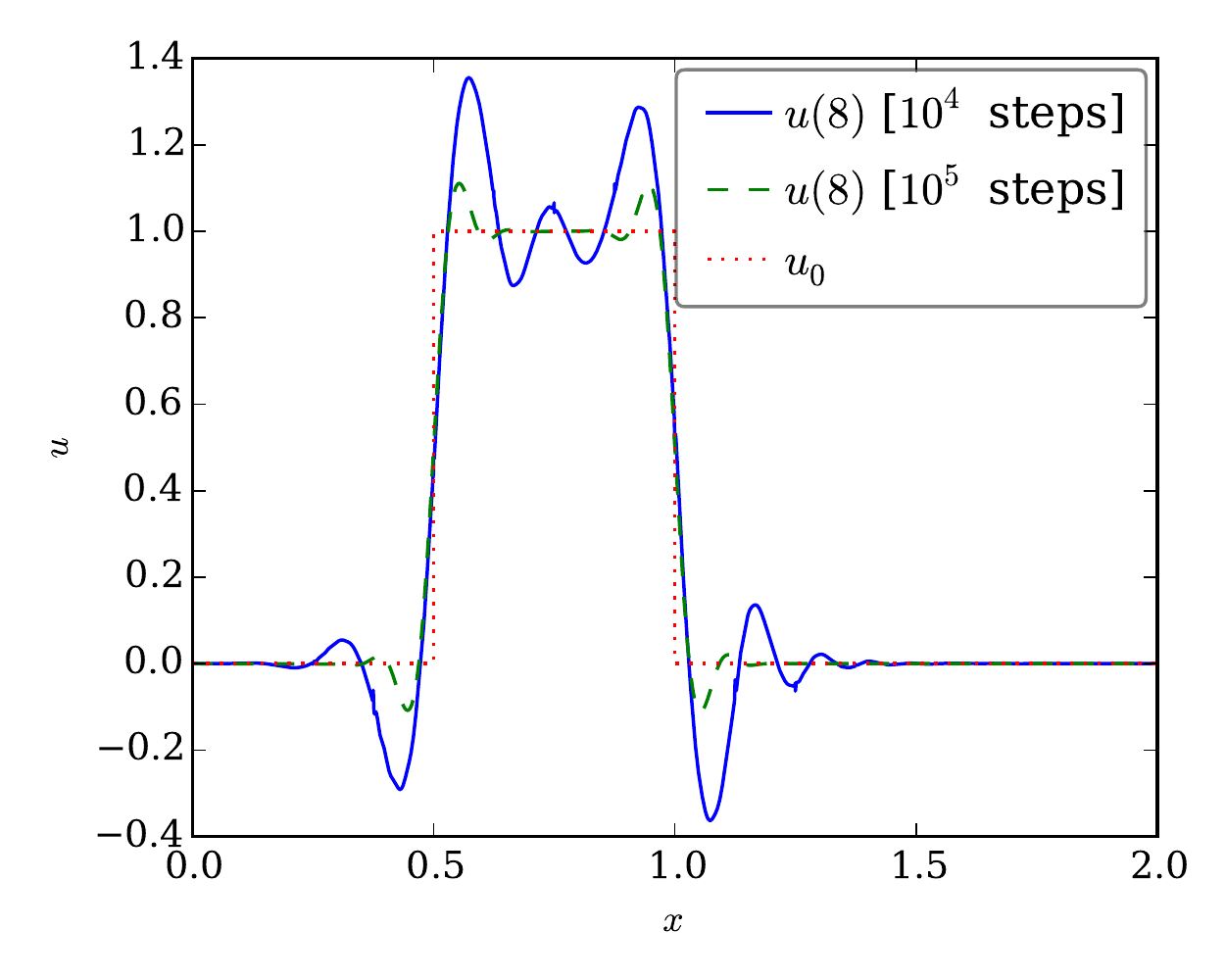}
    \caption{Solution computed using adaptive artificial dissipation.}
    \label{fig:linear_jump_ad_adaptive_solution}
  \end{subfigure}%
  \\
  \begin{subfigure}[b]{0.49\textwidth}
    \includegraphics[width=\textwidth]{%
      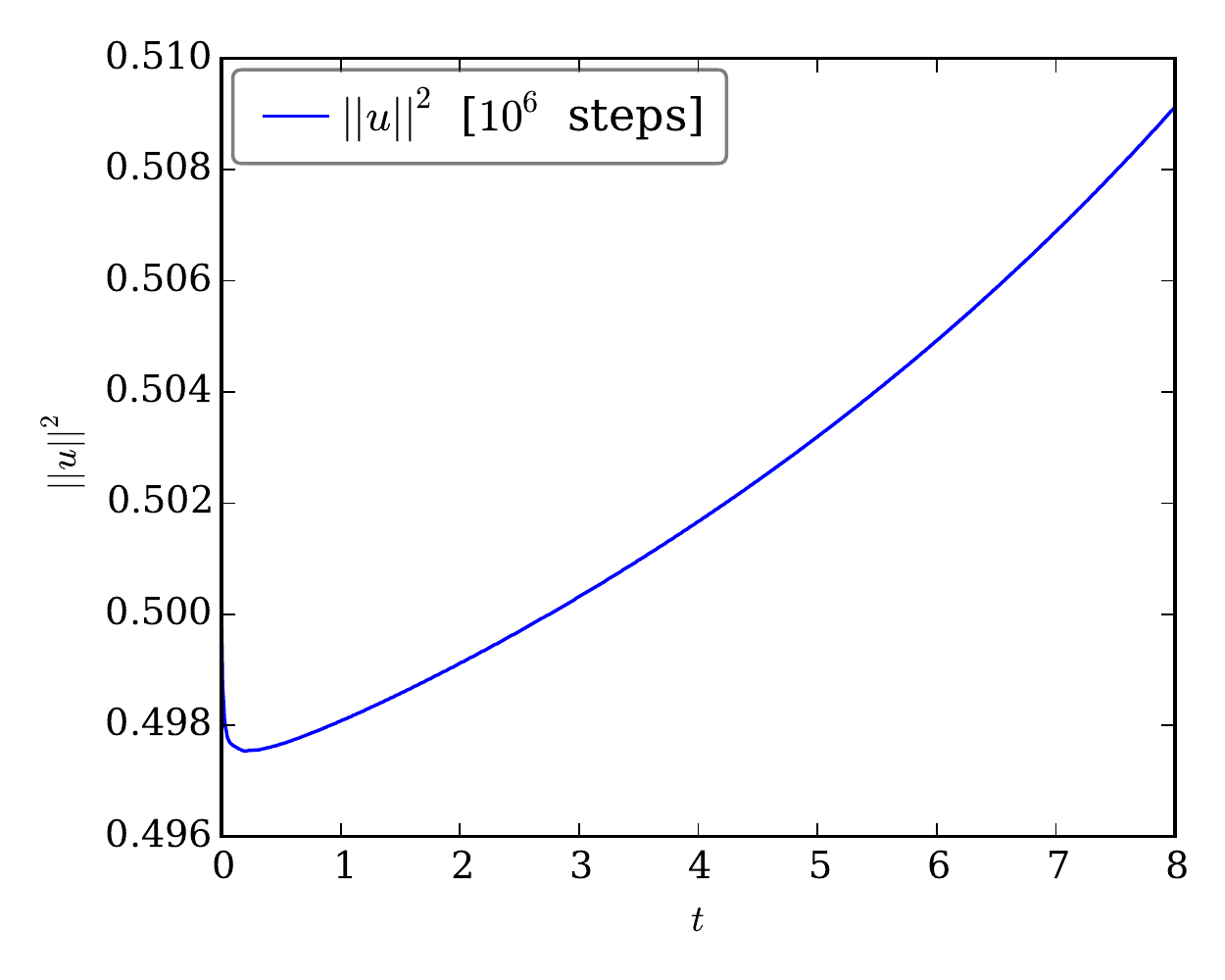}
    \caption{Energy of the solution computed without additional artificial dissipation.}
    \label{fig:linear_jump_no_filtering_energy}
  \end{subfigure}%
  ~
  \begin{subfigure}[b]{0.49\textwidth}
    \includegraphics[width=\textwidth]{%
      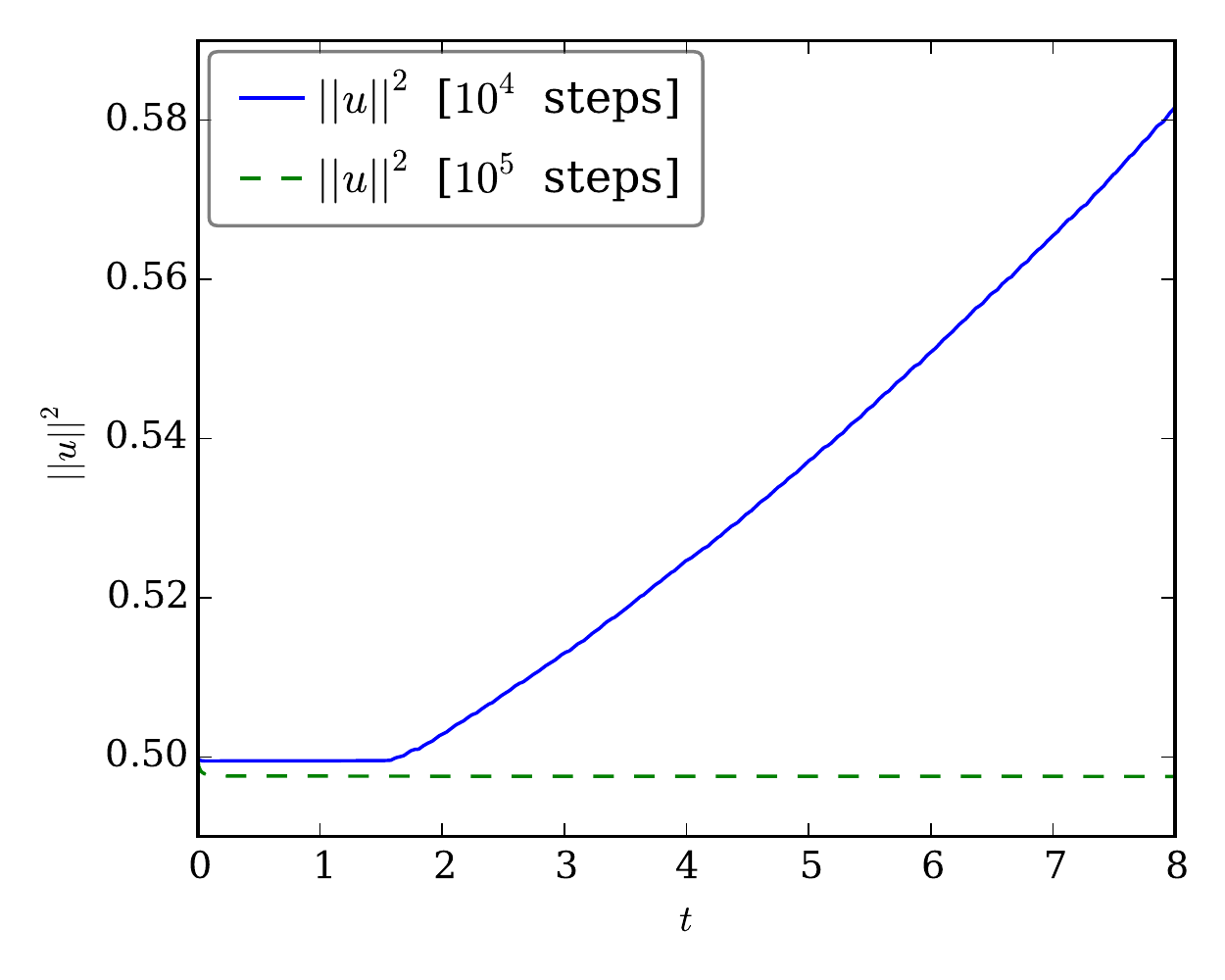}
    \caption{Energy of the solution computed using adaptive artificial dissipation.}
    \label{fig:linear_jump_ad_adaptive_energy}
  \end{subfigure}%
  \caption{Numerical results for linear advection using $N = 16$ elements with
           polynomials of degree $\leq p = 15$.
           On the left hand side, no additional artificial dissipation has been
           used, whereas adaptive spectral viscosity of orders $s=1$ has been
           used for the right hand side.
           }
  \label{fig:linear_jump}
\end{figure}

Using $10^6$ time steps, the energy of the solution computed without artificial
dissipation increases, as can be seen in Figure \ref{fig:linear_jump_no_filtering_energy}.
Additionally, Figure \ref{fig:linear_jump_no_filtering_solution} shows dominant
oscillations that have been developed. Contrary, adaptive artificial dissipation
stabilises the scheme. Therefore, $10^4$ time steps suffice to get a bounded
increase in the energy and some oscillations, see Figures
\ref{fig:linear_jump_ad_adaptive_energy} and \ref{fig:linear_jump_ad_adaptive_solution}.
With the same number of time steps, the computation without spectral viscosity
blows up. However, the artificial dissipation operator introduces a restriction
on the possible time steps. Thus, if the time step is too big, the desired
estimate on the strength $\epsilon$ fails. Setting $\epsilon$ to zero in this
case is the best possible solution, but the energy might increase.
This phenomenon is visible in Figure \ref{fig:linear_jump_ad_adaptive_energy}.
Increasing the number of time steps to $10^5$ yields the desired constant energy
and less oscillations in the solution plotted in Figure
\ref{fig:linear_jump_ad_adaptive_solution}.

\subsection{Burgers' equation}

Burgers' equation \eqref{eq:Burgers} with smooth initial condition
\begin{equation}
  \partial_t u + \partial_x \frac{u^2}{2} = 0
  ,\quad
  u(0,x) = u_0(x) = \sin \pi x + 0.01
\end{equation}
in the periodic domain $[0, 2]$ is used as a prototypical example of a nonlinear
conservation law yielding a discontinuous solution in finite time $t \in [0,3]$.
The stable semidiscretisation \eqref{eq:Burgers-semidisc} with $N = 16$ elements
representing polynomials of degree $\leq p = 16$ in nodal Gauß-Legendre bases
is used with the local Lax-Friedrichs flux $\fnum(u_-,u_+) = \frac{u_-^2 + u_+^2}{4}
- \frac{ \max \set{ \abs{u_-}, \abs{u_+} } }{2} (u_+ - u_-)$. The explicit Euler
method as time integrator uses $15 \cdot 10^3$ steps for the interval $[0, 3]$.

\begin{figure}[!htb]
  \centering
  \begin{subfigure}[b]{0.49\textwidth}
    \includegraphics[width=\textwidth]{%
      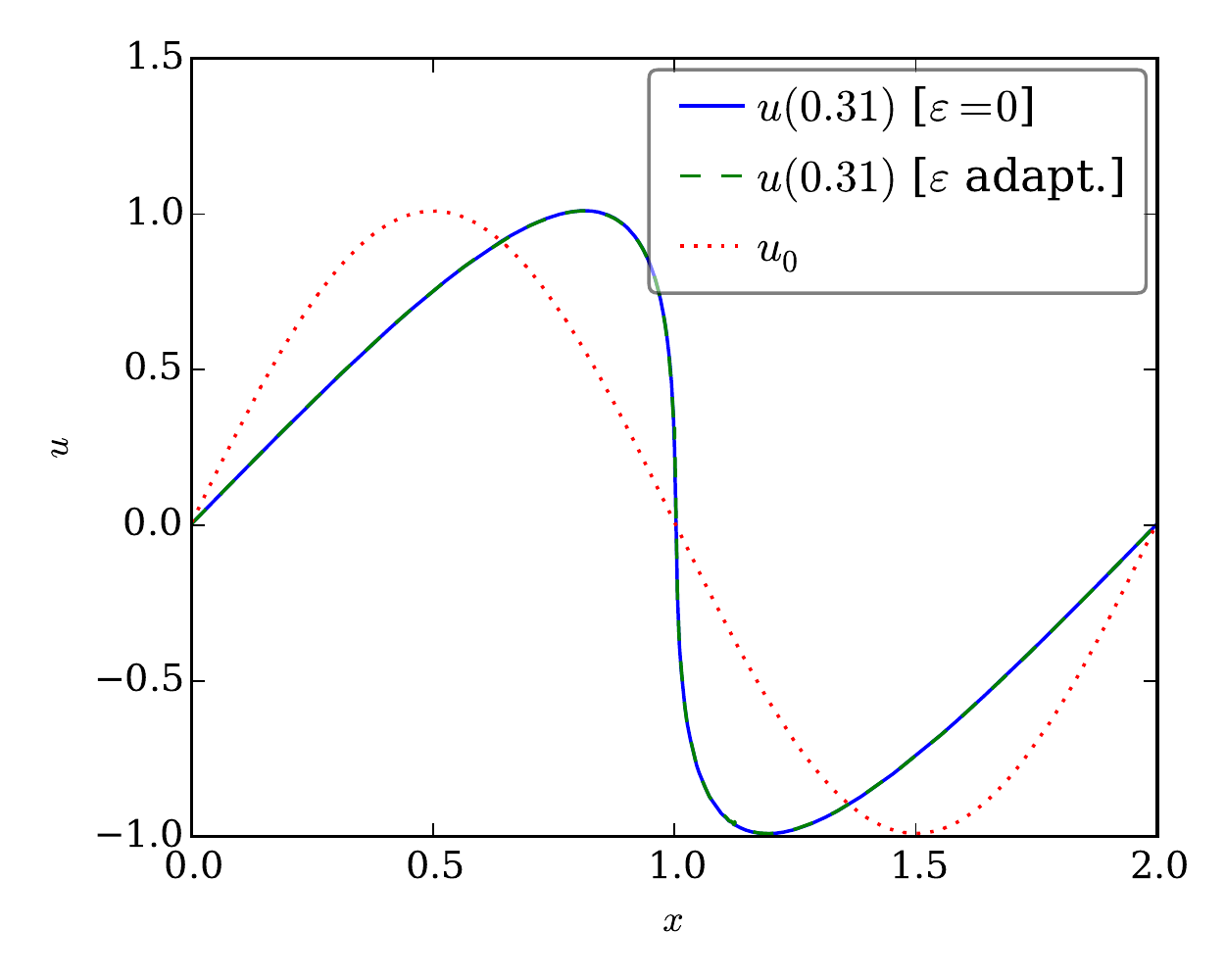}
    \caption{Solution at time $t = 0.31$.}
    \label{fig:burgers_sin_smooth_solution}
  \end{subfigure}%
  ~
  \begin{subfigure}[b]{0.49\textwidth}
    \includegraphics[width=\textwidth]{%
      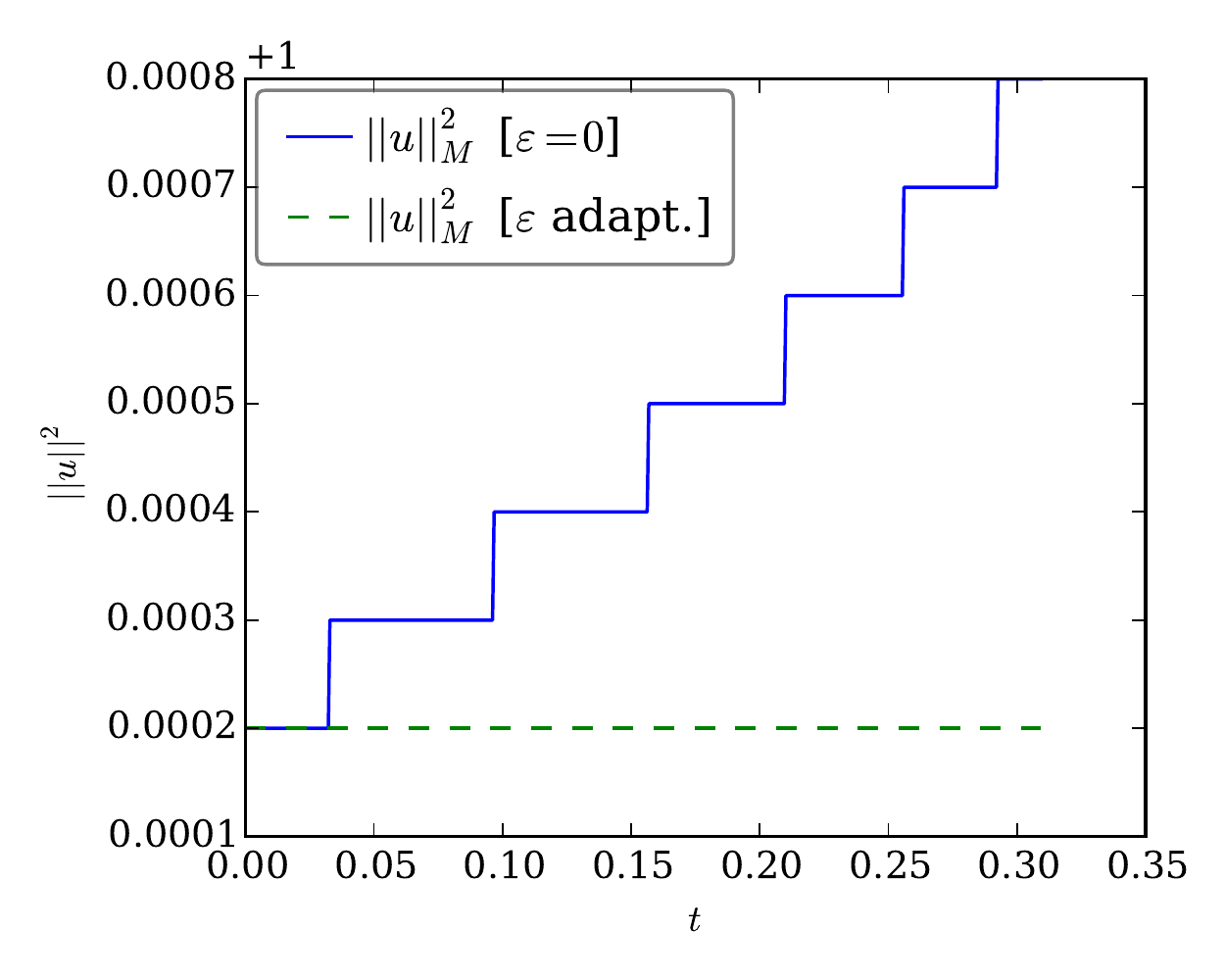}
    \caption{Energy in the time interval $[0, 0.31]$.}
    \label{fig:burgers_sin_smooth_energy}
  \end{subfigure}%
  \\
  \begin{subfigure}[b]{0.49\textwidth}
    \includegraphics[width=\textwidth]{%
      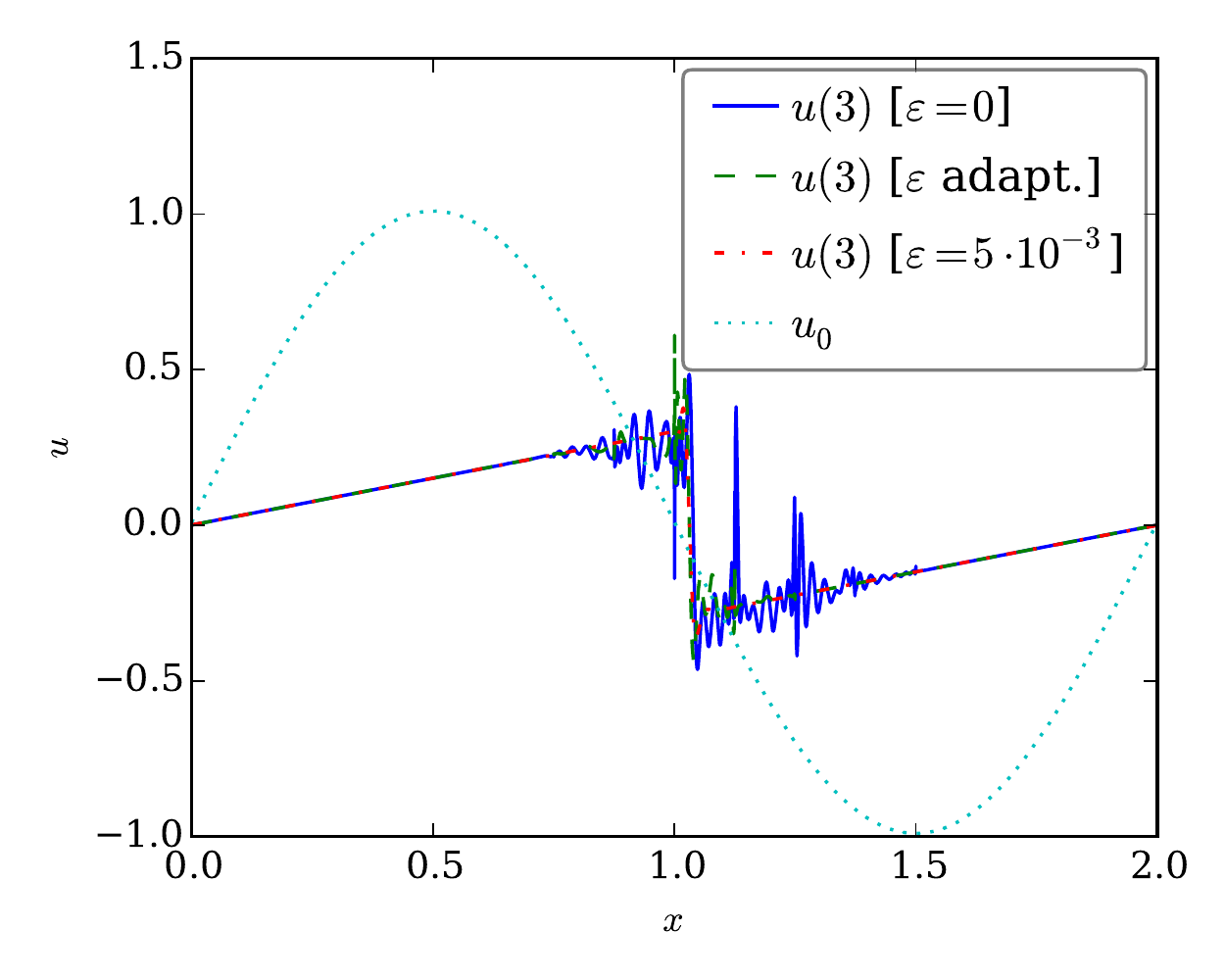}
    \caption{Solution at time $t = 3$.}
    \label{fig:burgers_sin_rough_solution}
  \end{subfigure}%
  ~
  \begin{subfigure}[b]{0.49\textwidth}
    \includegraphics[width=\textwidth]{%
      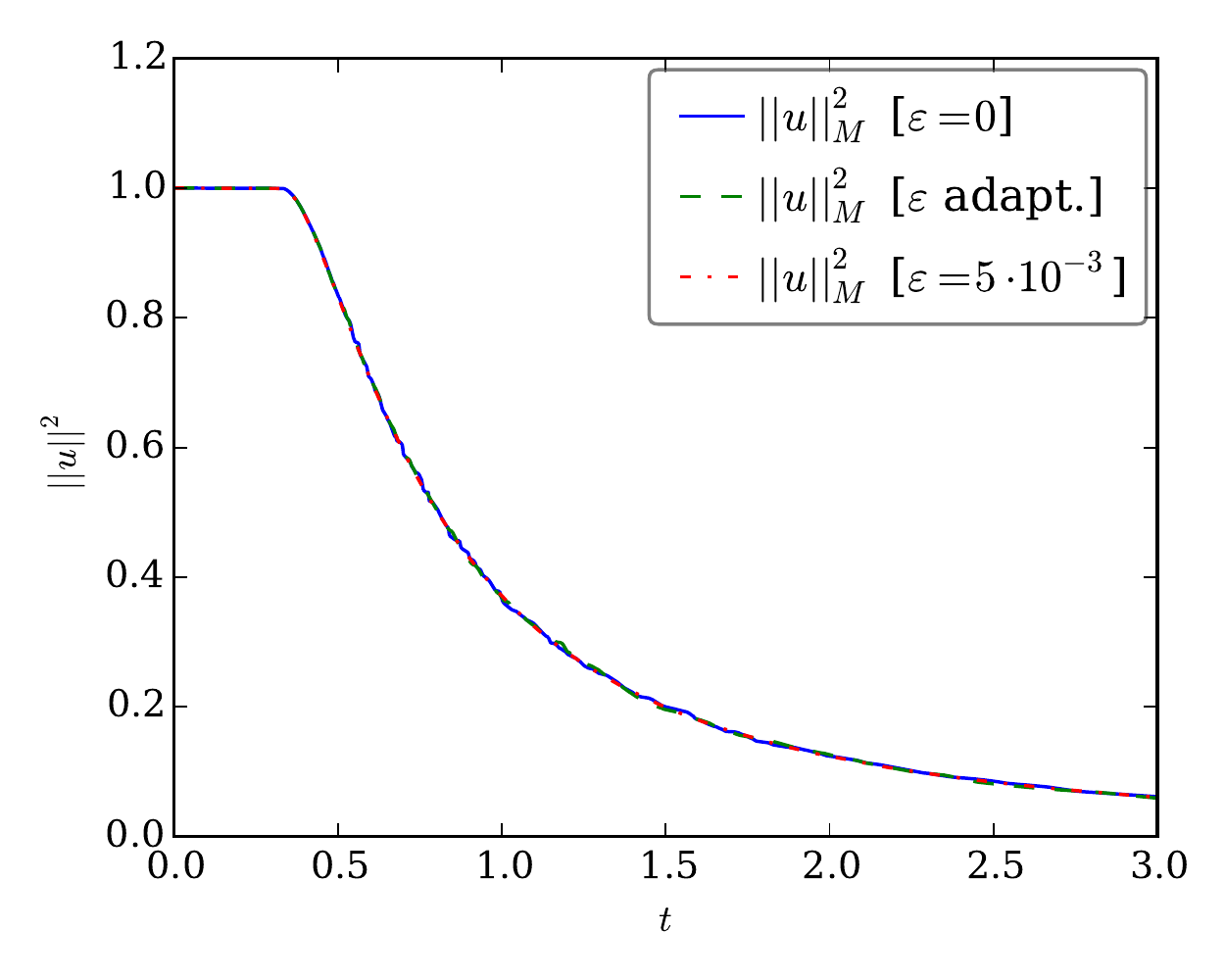}
    \caption{Energy in the time interval $[0, 3]$.}
    \label{fig:burgers_sin_rough_energy}
  \end{subfigure}%
  \caption{Numerical results for Burgers' equation using $N = 16$ elements with
           polynomials of degree $\leq p = 15$.
           The solutions and the energy are plotted on the left hand side and
           right hand side, respectively.
           }
  \label{fig:burgers_sin}
\end{figure}

At time $t = 0.31$, the solution in Figure \ref{fig:burgers_sin_smooth_solution}
computed with only $500$ time steps
is still smooth. However, the energy in Figure \ref{fig:burgers_sin_smooth_energy}
increases if no artificial dissipation is used. Contrary applying adaptive spectral
viscosity results in a constant energy.

At time $t = 3$, the solution in Figure \ref{fig:burgers_sin_rough_solution} has
developed a discontinuity resulting in oscillations around $x \approx 1$.
Although adaptive artificial dissipation damps these a bit, it does not removed
them. However, spectral viscosity of fixed strength $\epsilon = 5 \cdot 10^{-3}$
adds enough dissipation to remove them and yields a non-oscillatory result.
Nevertheless, all three choices of spectral viscosity yield nearly
visually indistinguishable results for the energy in Figure
\ref{fig:burgers_sin_rough_energy} due to the dissipative numerical flux.

  \section{Conclusions and further research}
\label{sec:conclusions}

In this work, artificial dissipation / spectral viscosity has been considered in
the general framework of CPR methods using SBP operators. A naive discretisation
of the viscosity operator does not yield the desired results, whereas the chosen
representation (after the application of summation-by-parts and cancellation of
undesired boundary terms) results in the estimates expected from the continuous
setting.

Additionally, a new adaptive strategy has been proposed in order to compute the
strength of the viscosity in a way to get a stable fully discrete scheme obtained
by an explicit Euler method. Thus, additional terms of order $(\Delta t)^2$ that
appear in the estimate of the energy growth in one time step are compensated.
However, this artificial dissipation is not enough to remove all oscillations,
especially the ones developing in nonlinear problems.

Numerical results for linear advection and Burgers' equation have been presented,
showing the advantages of the chosen approach as well as some limitations.
The application of artificial dissipation stabilises the scheme, but also
introduces additional restrictions on the time step. Therefore, removing these
by an operator splitting approach is desired and conducted in the second part
of this series \citex{in}{by} \citet{glaubitz2016enhancing}.

Another topic of further research is the investigation of different time
integration methods. While strong-stability preserving (SSP) schemes can be
written as convex combinations of explicit Euler steps and inherit therefore
the stability properties, other adaptive strategies might be advantageous in
this setting.

Moreover, extending the approach to other hyperbolic conservation laws will be
interesting.

  \section{Appendix: Legendre polynomials}

The Legendre polynomials can be represented by Rodrigues' formula
(equation 8.6.18 of \cite{abramowitz1972handbook})
\begin{equation}
  \phi_p(x) = \frac{1}{2^p \, p!} \od[p]{}{x} (x^2 - 1)^p
\end{equation}
and are orthogonal in $L_2[-1,1]$ with $\norm{\phi_p}^2 = 2 / (2p + 1)$. Their
boundary values are $\phi_p(1) = 1$ and $\phi_p(-1) = (-1)^p$. Due to Rodrigues'
formula, they are symmetric for even $p$ and antisymmetric for odd $p$.
Additionally, they obey
\begin{equation}
\label{eq:Legendre-derivative}
\begin{aligned}
  \phi_{p+1}'(x)
    &= \frac{1}{2^{p+1} \, (p+1)!} \od[p+2]{}{x} (x^2 - 1)^{p+1}
  \\&= \frac{1}{2^{p+1} \, (p+1)!} \od[p+1]{}{x} \left[ 2 (p+1) x (x^2 - 1)^{p} \right]
  \\&= \frac{1}{2^{p} \, p!} \od[p]{}{x} \od{}{x} \left[ x (x^2 - 1)^{p} \right]
  \\&= \frac{1}{2^{p} \, p!} \od[p]{}{x} \left[(x^2 - 1)^{p} + 2 p x^2 (x^2 - 1)^{p-1}\right]
  \\&= \frac{1}{2^{p} \, p!} \od[p]{}{x} \left[(2p+1) (x^2 - 1)^{p} + 2 p (x^2 - 1)^{p-1}\right]
  \\&= (2p+1) \frac{1}{2^{p} \, p!} \od[p]{}{x} (x^2 - 1)^{p}
       + \frac{1}{2^{p-1} \, (p-1)!} \od[p]{}{x} (x^2 - 1)^{p-1}
  \\&= (2p+1) \phi_p(x) + \phi_{p-1}'(x).
\end{aligned}
\end{equation}
The first three Legendre polynomials are $\phi_0(x) = 1$, $\phi_1(x) = x$,
$\phi_2(x) = (3 x^2 - 1) / 2$.

In order to compute the $L_2$ projection of $(1-x^2) \phi_p(x)$ on the space of
polynomials of degree $\leq p$, equation 8.5.4 of \citet{abramowitz1972handbook}
can be used
\begin{equation}
  (x^2 - 1) \partial_x \phi_p(x) = p x \phi_p(x) - p \phi_{p-1}(x).
\end{equation}
Inserting equation 8.5.3 of \citet{abramowitz1972handbook}
\begin{equation}
  (p+1) \phi_{p+1}(x) = (2p+1) x \phi_p(x) - p \phi_{p-1}(x)
\end{equation}
results in
\begin{equation}
\label{eq:RHS-Legendre-p}
\begin{aligned}
  (x^2 - 1) \partial_x \phi_p(x)
  &= \frac{p (p+1)}{2p+1} \phi_{p+1} + \frac{p^2}{2p+1} \phi_{p-1}(x) - p \phi_{p-1}(x)
  \\&= \frac{p (p+1)}{2p+1} \phi_{p+1}(x) - \frac{p (p+1)}{2p+1} \phi_{p-1}(x).
\end{aligned}
\end{equation}

The Lobatto-Legendre quadrature includes both boundary nodes and is exact
for polynomials of degree $\leq 2p-1$. The norm of $\phi_p$, evaluated by
Lobatto-Legendre quadrature is (equation (1.136) of \citet{kopriva2009implementing})
\begin{equation}
\label{eq:norm-lobatto-phi_p}
  \norm{\vec{\phi}_p}_{M}^2
  = \vec{\phi}_p^T \mat{M} \vec{\phi}_p
  = \frac{2}{p}.
\end{equation}
In order to compute $\vec{\phi}_{p-1}^T \mat{M} \vec{\phi}_{p+1}$ via Lobatto-Legendre
quadrature, the product can be expanded as a linear combination of Legendre polynomials
\begin{equation}
  \phi_{p-1} \phi_{p+1} = \sum_{n = 0}^{2p} \beta_n \phi_n.
\end{equation}
Since the Legendre polynomials are orthogonal, $\beta_0 = 0$. As used \citex{in}{by}
\citet{vincent2011newclass}, the leading coefficient of the Legendre polynomial
$\phi_n$ of degree $n$ is
\begin{equation}
\label{eq:leading-coefficient-legendre}
  a_n = \frac{(2n)!}{2^n (n!)^2}.
\end{equation}
Thus,
\begin{equation}
\begin{aligned}
  \beta_{2p} &= \frac{a_{p-1} a_{p+1}}{a_{2p}}
  = \frac{(2p-2)!}{2^{p-1} ((p-1)!)^2}
    \frac{(2p+2)!}{2^{p+1} ((p+1)!)^2}
    \frac{2^{2p} ((2p)!)^2}{(4p)!}
  \\&
  = \frac{ (2p-2)! \, (2p+2)! \, ((2p)!)^2 }{ ((p-1)!)^2 \, ((p+1)!)^2 \, (4p)!}.
\end{aligned}
\end{equation}
Denoting the approximation of $\int \cdot$ by Lobatto-Legendre quadrature as
$\int_L \cdot$, $\int_L \phi_{2p}$ has to be computed. To use equation
\eqref{eq:norm-lobatto-phi_p}, $\phi_p^2$ is expanded similar to $\phi_{p-1} \phi_{p+1}$
\begin{equation}
  \phi_p^2 = \sum_{n = 0}^{2p} \gamma_n \phi_n.
\end{equation}
Since Legendre polynomials are orthogonal, 
\begin{equation}
  \gamma_0
  = \norm{\phi_p}^2 / \int \phi_0
  = \frac{1}{2p+1}.
\end{equation}
Similar to $\beta_{2p}$, $\gamma_{2p}$ can be written as
\begin{equation}
  \gamma_{2p} = \frac{a_p^2}{a_{2p}}
  = \frac{\left[(2p)!\right]^2}{2^{2p} (p!)^4}
    \frac{2^{2p} ((2p)!)^2}{(4p)!}
  = \frac{ \left[(2p)!\right]^4}{ (p!)^4 \, (4p)!}.
\end{equation}
Using $\int_L \phi_p^2 = 2 / p$ from equation \eqref{eq:norm-lobatto-phi_p}
and linearity of the quadrature yields
\begin{equation}
  \frac{2}{p}
  = \int_L \phi_p^2
  = \sum_{n = 0}^{2p} \gamma_n \int \phi_n
  = \gamma_0 \int_L \phi_0 + \gamma_{2p} \int_L \phi_{2p}
  = 2 \gamma_0 + \gamma_{2p} \int_L \phi_{2p},
\end{equation}
since the quadrature is exact (and thus zero) for $\phi_n, n \in \set{1, \dots, 2p-1}$.
Therefore,
\begin{equation}
  \int_L \phi_{2p}
  = \left( \frac{2}{p} - 2 \gamma_0 \right) \gamma_{2p}^{-1}
  = \left( \frac{2}{p} - \frac{2}{2p+1} \right) \gamma_{2p}^{-1}
  = 2 \frac{p+1}{p (2p+1)} \frac{ (p!)^4 \, (4p)!}{ \left[(2p)!\right]^4}.
\end{equation}
Finally, using $\beta_0 = 0$ and $\int_L \phi_n = 0$ for $n \in \set{1, \dots, 2p-1}$,
\begin{equation}
\label{eq:lobatto-phi_p-1-phi_p+1}
\begin{aligned}
  \vec{\phi}_{p-1}^T \mat{M} \vec{\phi}_{p+1}
  &= \int_L \phi_{p-1} \phi_{p+1}
  = \sum_{n = 0}^{2p} \beta_n \int_L \phi_n
  = \beta_{2p} \int_L \phi_{2p}
  \\&
  = \frac{ (2p-2)! \, (2p+2)! \, ((2p)!)^2 }{ ((p-1)!)^2 \, ((p+1)!)^2 \, (4p)!}
    \cdot 2 \frac{p+1}{p (2p+1)} \frac{ (p!)^4 \, (4p)!}{ \left[(2p)!\right]^4}
  \\&
  = \frac{2}{p}
    \frac{p+1}{2p+1}
    \frac{ (p!)^2 }{ ((p-1)!)^2 }
    \frac{ (p!)^2 }{ ((p+1)!)^2 }
    \frac{ (2p-2)! }{ (2p)! }
    \frac{ (2p+2)! }{ (2p)! }
  \\&
  = \frac{2}{p}
    \frac{p+1}{2p+1}
    \frac{ p^2 }{ 1 }
    \frac{ 1 }{ (p+1)^2 }
    \frac{ 1 }{ 2 p (2p-1) }
    \frac{ (2p+2) (2p+1) }{ 1 }
  \\&
  = \frac{2}{p}
    \frac{ p }{ 2p-1 }.
\end{aligned}
\end{equation}

  \printbibliography

\end{document}